\documentclass[aos]{imsart}
\usepackage{xcite}

\makeatletter
\newcommand*{\addFileDependency}[1]{
  \typeout{(#1)}
  \@addtofilelist{#1}
  \IfFileExists{#1}{}{\typeout{No file #1.}}
}
\makeatother

\RequirePackage{amsthm,amsmath,amsfonts,amssymb}
\RequirePackage{dsfont}
\RequirePackage{color}

\RequirePackage[colorlinks,citecolor=blue]{hyperref}
\RequirePackage[authoryear]{natbib} %% uncomment this for author-year bibliography
\RequirePackage{enumerate}
\RequirePackage{enumitem}
\RequirePackage{graphicx}% uncomment this for including figures
\allowdisplaybreaks[4]
\startlocaldefs
\theoremstyle{plain}
\newtheorem{theorem}{Theorem}
\newtheorem{lemma}{Lemma}

\newtheorem{corollary}{Corollary}
\theoremstyle{remark}

\newtheorem{assum}{Assumption}
\newcommand{\var}{\mathbb{V}\mathrm{ar}}
\newcommand{\diff}{\mathrm{d}}
\newcommand{\E}{\mathbb{E}}
\newcommand{\K}{\mathrm{K}}
\newcommand{\rhsnorm}{\|_{\mathrm{HS}} }
\newcommand{\cov}{\operatorname{Cov}}
\newcommand{\Lcal}{\mathcal{L}}

\endlocaldefs

\begin{document}

\begin{frontmatter}
\title{Theory of functional principal component analysis \\ for discretely observed data}
%\title{A sample article title with some additional note\thanksref{T1}}
\runtitle{Theory of FPCA for discretely observed data}
%\thankstext{T1}{A sample of additional note to the title.}

\begin{aug}
%%%%%%%%%%%%%%%%%%%%%%%%%%%%%%%%%%%%%%%%%%%%%%
%%Only one address is permitted per author. %%
%%Only division, organization and e-mail is %%
%%included in the address.                  %%
%%Additional information can be included in %%
%%the Acknowledgments section if necessary. %%
%%%%%%%%%%%%%%%%%%%%%%%%%%%%%%%%%%%%%%%%%%%%%%
\author[A]{\fnms{Hang} \snm{Zhou}\ead[label=e1]{hgzhou@ucdavis.edu}},
\author[B]{\fnms{Dongyi} \snm{Wei}\ead[label=e2]{jnwdyi@pku.edu.cn}}
\and
\author[C]{\fnms{Fang} \snm{Yao}\ead[label=e3]{fyao@math.pku.edu.cn} }\footnote{Fang Yao is the corresponding author. }
%%%%%%%%%%%%%%%%%%%%%%%%%%%%%%%%%%%%%%%%%%%%%%
%% Addresses                                %%
%%%%%%%%%%%%%%%%%%%%%%%%%%%%%%%%%%%%%%%%%%%%%%
\address[A]{Department of Statistics, University of California at Davis, \printead{e1}}
\address[B]{School of Mathematical Sciences, Peking University, \printead{e2}}
\address[C]{School of Mathematical Sciences and Center for Statistical Science, Peking University, \printead{e3}}

\end{aug}

\begin{abstract}
Functional data analysis is an important research field in statistics which treats  data as random functions drawn from some infinite-dimensional functional space, and functional principal component analysis (FPCA) based on eigen-decomposition plays a central role for data reduction and representation. After nearly three decades of research, there remains a key problem unsolved, namely,  the perturbation analysis of covariance operator for diverging number of eigencomponents obtained from noisy and discretely observed data. This is fundamental for studying models and methods based on FPCA, while there has not been substantial progress since \cite{hall2006propertieshans}'s result for a fixed number of eigenfunction estimates. In this work, we aim to establish a unified theory for this problem, obtaining upper bounds for eigenfunctions with diverging indices in both the $\mathcal{L}^2$ and supremum norms, and deriving the asymptotic distributions of eigenvalues for a wide range of sampling schemes. Our results provide insight into the phenomenon when the $\mathcal{L}^{2}$ bound of eigenfunction estimates with diverging indices is minimax optimal as if the curves are fully observed, and reveal the transition of convergence rates from nonparametric to parametric regimes in connection to sparse or dense sampling.  We also develop a double truncation technique to handle the uniform convergence  of estimated covariance and eigenfunctions. The technical arguments in this work are useful for handling the perturbation series with noisy and discretely observed functional data and can be applied in models or those involving inverse problems based on FPCA as regularization, such as functional linear regression.
\end{abstract}



\end{frontmatter}
%%%%%%%%%%%%%%%%%%%%%%%%%%%%%%%%%%%%%%%%%%%%%%
%% Please use \tableofcontents for articles %%
%% with 50 pages and more                   %%
%%%%%%%%%%%%%%%%%%%%%%%%%%%%%%%%%%%%%%%%%%%%%%
%\tableofcontents

%%%%%%%%%%%%%%%%%%%%%%%%%%%%%%%%%%%%%%%%%%%%%%
%%%% Main text entry area:

\section{Introduction}\label{sec:intro}

Modern data collection technologies have rapidly evolved, leading to the widespread emergence of functional data that have been extensively studied over the past few decades. Generally, functional data are considered stochastic processes that satisfy certain smoothness conditions or {realizations of random elements valued in Hilbert space.} These two perspectives highlight the essential natures of functional data, {namely, their smoothness and infinite dimensionality, which distinguish} them from high-dimensional and vector-valued data. For a comprehensive treatment of functional data, we recommend the monographs by \cite{ramsay2005functional}, \cite{ferraty2006nonparametric}, \cite{horvath2012inference}, {and} \cite{hsing2015theoretical}, among others.

Although functional data provide information over a continuum, which is often time or spatial location, in reality, data are collected or observed discretely  with measurement errors. For instance, we usually use $n$ to denote the sample size, which is the number of subjects corresponding to random functions, and $N_i$ to denote the number of observations for the $i$th subject. Thanks to the smooth nature of functional data, {having a large number of observations per subject is more of a blessing than a curse,} in contrast to high-dimensional data \citep{hall2006propertieshans, zhang2016sparse}. There is extensive literature on nonparametric methods that address the smoothness of functional data, including kernel or local polynomial methods \citep{yao2005jasa, hall2006propertieshans, zhang2016sparse}, and various types of spline methods \citep{rice2001nonparametric, yao2006penalized, paul2009consistency, cai2011optimal}.

When employing a smoothing method, there are two typical strategies to be considered. If the observed {time} points per subject are relatively dense, it is recommended to pre-smooth each curve before further analysis, as suggested by \cite{ramsay2005functional} and \cite{zhang2007statistical}. However, if the sampling scheme is rather sparse, {it is preferred to pool observations together from all subjects \citep{yao2005jasa}.} The choice of individual versus pooled information affects the convergence rates and phase transitions in estimating population quantities, such as mean and covariance functions. When $N_i\gtrsim O(n^{5/4})$ and the tuning parameter is optimally chosen per subject, the estimated mean and covariance functions based on the reconstructed curves {through} pre-smoothing are $\sqrt{n}$-consistent, the so-called optimal parametric rate. On the other hand, by borrowing information from all subjects, the pooling method only requires $N_i\gtrsim O(n^{1/4})$ for mean and covariance estimation to reach optimal \citep{cai2010nonparametric, cai2011optimal, zhang2016sparse}, providing theoretical insight into the advantage of the pooling strategy.

However, estimating {the} mean and covariance functions does not account for the infinite dimensionality of functional data. Due to the non-invertibility of covariance operators for functional random objects, regularization is required in models that involve  {inverse issues} with functional covariates, such as the functional linear model \citep{yao2005aos, hall2007methodology, yuan2010reproducing}, functional generalized linear model \citep{muller2005generalized, dou2012estimation}, and functional Cox model \citep{qu2016optimal}. Truncation of the leading functional principal components (FPC) is a well-developed approach to addressing {this} inverse issue \citep{hall2007methodology, dou2012estimation}. In order to suppress the model bias, the number of principal components used in truncation should grow slowly with sample size. As a result, the convergence rate of the estimated eigenfunctions with diverging indices {becomes} a fundamental issue, which is not only important in its own right {but  also crucial for most models and methods involving functional principal components regularization.} 

For fully observed functional data, \cite{hall2007methodology} obtained the optimal convergence rate  $j^2/n$ for the $j$th eigenfunction, which served as a cornerstone in establishing the optimal convergence {rate} in functional linear model \citep{hall2007methodology} and functional generalized linear model \citep{dou2012estimation}. In the discretely observed case, stochastic bounds for a fixed number of eigenfunctions have been obtained by different methods. Using {a} local linear smoother, \cite{hall2006propertieshans} showed that the $\Lcal^2$ rate of a fixed eigenfunction for finite $N_i$ is $O_{P}(n^{-4/5})$. Under the reproducing kernel Hilbert space framework, \cite{cai2010nonparametric} claimed that eigenfunctions with fixed indices admit the same convergence rate as the covariance function, which is $O_{P}((n/\log n)^{-4/5})$. It is important to note that, although both results are one-dimensional nonparametric rates ({differing} at most by a factor of $(\log n)^{4/5}$), the methodologies and techniques used are completely disparate, and a detailed discussion can be found in Section \ref{sec:preli}. Additionally, \cite{paul2009consistency} proposed a reduced rank model and studied its asymptotic properties under a particular setting. 
{In \cite{zhou2023functional}, the authors studied the convergence rate for the functional linear model and obtained an improved bound for the eigenfunctions with diverging indices. However, this rate will not reach the optimal rate of $j^{2}/n$ for any sampling rate $N_i$. As explained  in Section \ref{sec:preli}, while some bounds can be obtained for eigenfunctions with diverging indices, attaining an optimal bound presents a substantially greater challenge. Lack of such an optimal bound for eigenfunctions poses considerable challenge in analyzing the standard and efficient plug-in estimator in functional linear model \citep{hall2007methodology}. Consequently, \cite{zhou2023functional} resorted to a complex sample-splitting strategy, which results in lower estimation efficiency.}
To the best of our knowledge, there has been no progress in obtaining the {optimal} convergence rate of eigenfunctions with diverging indices when the data are discretely observed with noise contamination.

The distinction between estimating {a diverging number and a fixed number of eigenfunctions} is rooted in the infinite-dimensional nature of functional data. Analyzing eigenfunctions with diverging indices {presents challenges} due to the decaying eigenvalues. For fully observed data, the cross-sectional sample covariance based on {the} true functions facilitates the {application} of perturbation results, as shown in prior work \citep{hall2007methodology, dou2012estimation}. This approach simplifies each term in the perturbation series to the principal component scores. However, when the trajectories are observed at discrete time points, this virtue no longer exists, {leading to a summability issue arising from the estimation bias and decaying eigenvalues.} {This  renders existing techniques invalid and remains an unresolved problem; see Section \ref{sec:preli} for further elaboration.}

{This paper addresses this significant yet challenging task of estimating an increasing number of eigenfunctions from discretely observed functional data, and presents a unified theory. The contributions of this paper are at least threefold. First, we establish an $\mathcal{L}^2$ bound for the eigenfunctions and the asymptotic normality of the eigenvalues with increasing indices, reflecting a transition from nonparametric to parametric regimes and encompassing a wide range from sparse to dense sampling. We show that when $N_i$ reaches a magnitude of $n^{1/4+\delta}$, where $\delta$ depends on the  smoothness parameters of the underlying curves, the convergence rate becomes optimal as if the curves are fully observed.  Second, we introduce a novel double truncation method that yields uniform convergence across the time domain, surmounting theoretical barriers in the existing literature. Through this approach, uniform convergence rates for the covariance and eigenfunctions are achieved under mild conditions across various sampling schemes. Notably, this includes the uniform convergence of eigenfunctions with increasing indices, which is new even in scenarios where data are fully observed. Third, we provide a new technical route for addressing the perturbation series of the functional covariance operator, bridging the gap between the ``ideal'' fully observed scenario and the noisy, discrete ``real-world'' context. These advanced techniques pave the way for their application in downstream FPCA-related analyses, and the achieved optimal rate of eigenpairs facilitates the extension of existing theoretical results for ``fully'' observed functional data to discreetly observed case.}

The rest of the paper is organized as follows. In Section \ref{sec:preli}, we give a synopsis of  covariance and eigencomponents estimation in functional data. We present the $\mathcal{L}^{2}$ convergence of eigenfunctions in Section \ref{sec:eigthm}, and discuss the uniform convergence problem of functional data in Section \ref{sec:uni}. Asymptotic {normality} of eigenvalues is presented in Section \ref{sec:asy}. {Section \ref{sec:simu} provides an illustration of the phase transition phenomenom in eigenfunctions with synthetic data.} The proofs of Theorem \ref{thm-eig} can be found in Appendix, while the proofs of other theorems and lemmas are collected in the Supplementary Material.

In what follows, we denote by $A_{n}=O_{P}(B_{n})$ the relation $\mathbb{P}(A_{n}\leq MB_{n})\geq1-\epsilon$, and by $A_{n}=o_{P}(B_{n})$ the relation $\mathbb{P}(A_{n}\leq \epsilon B_{n})\rightarrow0$ as $n\rightarrow\infty$, for each $\epsilon>0$ and a positive constant $M$. A non-random sequence $a_{n}$ is said to be $O(1)$ if it is bounded. For any non-random sequence $b_{n}$, we say $b_{n}=O(a_{n})$ if $b_{n}/a_{n}=O(1)$, and $b_{n}=o(a_{n})$ if $b_{n}/a_{n}\rightarrow 0$.  The notation $a_n\lesssim b_n$ indicates $a_n\leq C b_n$ for sufficiently large $n$ and postive constant $C$, and the relation $\gtrsim$ is defined similarly. We write $a_n\asymp b_n$ if $a_n\lesssim b_n$ and $b_n\lesssim a_n$. For $a\in\mathbb{R}$, $\lfloor a \rfloor $ denotes the largest integer less than or equal to $a$. For a function $f \in \Lcal^{2}[0,1]$, where $\Lcal^{2}[0,1]$ denotes the space of square-integrable functions on $[0,1]$, $\|f\|^2$ denotes $\int_{[0,1]}f(s)^2\diff s $, and $\|f\|_{\infty}$ denotes $\sup_{s\in[0,1]}|f(s)|$. For a function $A(s,t)\in\Lcal^2[0,1]^{2}$, define $\|A\rhsnorm^{2}=\iint_{[0,1]^{2}}A(s,t)^2\diff s\diff t $ and $\|A\|_{(j)}^{2}=\int_{[0,1]}\{\int_{[0,1]}A(s,t)\phi_{j}(s)\diff s \}^2\diff t $, where $\{\phi_{j}\}_{j=1}^{\infty} $ are the  eigenfunctions of interest.  We write $\int pq$ and $\int Apq$ for $\int p(u)q(u)\diff u $ and $\iint A(u,v)p(u)q(v)\diff u\diff v $ occasionally for brevity.

\section{Eigen-estimation for discretely observed functional data}\label{sec:preli}

Let $X(t)$ be a square integrable stochastic process on $[0,1]$, and let $X_{i}(t)$ be independent and identically distributed copies of $X(t)$. The mean and covariance functions of $X(t)$ are denoted by $\mu(t)=\mathbb{E}[{X(t)}]$ and $C(s,t)=\E [{ \{X(s)-\mu(s)\}\{X(t)-\mu(t)\}}]$, respectively. According to Mercer's Theorem \citep{indritz1963methods}, $C(s,t)$ has the spectral decomposition
\begin{equation}\label{eq:mercer}
	C(s,t)=\sum_{k=1}^{\infty}\lambda_{k}\phi_{k}(s)\phi_{k}(t),
\end{equation}  
where $\lambda_{1}>\lambda_{2}>\ldots>0$ are eigenvalues and $\{\phi_{j}\}_{j=1}^{\infty}$ are the corresponding eigenfunctions, which  form a complete orthonormal system  on $\Lcal^{2}[0,1]$. For each $i$, the process $X_{i}$ admits the so-called Karhunen-Lo\`eve expansion
\begin{equation}\label{eq:KL}
	X_{i}(t)=\mu(t)+\sum_{j=1}^{\infty}\xi_{ik}\phi_{k}(t),
\end{equation} 
where $\xi_{ik}=\int_{0}^{1}\{X_{i}(t)-\mu(t)\}\phi_{k}(t)\diff t$ are functional principal component scores with  zero mean and  variances $\lambda_{k}$.

However, in practice, it is only an idealization to have each $X_{i}(t)$ for all $t\in[0,1]$ to simplify theoretical analysis. Measurements are typically taken at $N_{i}$ discrete time points with noise contamination. Specifically, the actual observations for each $X_{i}$ are given by
\begin{equation}\label{eq:Xmodel}
\{(t_{ij},X_{ij}) | X_{ij}=X_{i}(t_{ij})+\varepsilon_{ij}, j=1,\cdots,N_{i} \},
\end{equation}
where $\varepsilon_{ij}$ are random copies of $\varepsilon$, with $\mathbb{E}(\varepsilon)=0$ and $\text{Var}(\varepsilon )=\sigma_{X}^{2}$. {We further assume the measurements errors $\{\varepsilon_{ij}\}_{i,j}$ are independent of $X_i$.}

Local linear regression is a popular smoothing technique in functional data analysis due to its attractive theoretical properties \citep{yao2005jasa, hall2006propertieshans, li2010uniform, zhang2016sparse}. The primary goal of this paper is to develop a unified theory for estimating a larger number of eigenfunctions from discretely observed functional data. To maintain focus and avoid distractions, we assume that the mean function $\mu(t)$ is known, and set $\mu(t)=0$ without loss of generality. The scenario involving an unknown mean function is discussed later in Section \ref{sec:eigthm}. We denote by $ \delta_{ijl}=X_{ij}X_{il}$ the raw covariance, and define $v_{i}=\{nN_{i}(N_{i}-1)\}^{-1}$. The local linear estimator of the covariance function is given by $\hat{C}(s,t)=\hat{\beta}_{0}$,
\begin{equation}\label{eq:ll}
	\begin{aligned}
\left(\hat{\beta}_0, \hat{\beta}_1, \hat{\beta}_2\right)= & \underset{\beta_0, \beta_1, \beta_2}{\operatorname{argmin}} \sum_{i=1}^n v_i \sum_{1 \leq l_1 \neq l_2 \leq N_i}\left\{\delta_{i j l}-\beta_0-\beta_1\left(t_{i l_1}-s\right)-\beta_2\left(t_{i l_2}-t\right)\right\}^2 \\
& \quad\quad\quad \times\frac{1}{h} \K\left(\frac{t_{i l_1}-s}{h} \right)  \frac{1}{h} \K\left(\frac{t_{i l_2}-t}{h} \right) ,
\end{aligned}
\end{equation}
where $\mathrm{K}$ is a symmetric, Lipschitz continuous density kernel on $[-1,1]$ and $h$ is the tuning parameter. The estimated covariance function $\hat{C}(s,t)$ can be expressed as an empirical version of the spectral decomposition in \eqref{eq:mercer}, i.e.,
\begin{equation}\label{eq:emercer}
	\hat{C}(s,t)=\sum_{k=1}^{\infty}\hat{\lambda}_{k}\hat{\phi}_{k}(s)\hat{\phi}_{k}(t),
\end{equation} 
where $\hat{\lambda}_{k}$ and $\hat{\phi}_{k}$ are estimators for $\lambda_{k}$ and $\phi_{k}$, respectively.  We assume that $\langle\hat{\phi}_{k},\phi_{k} \rangle\geq 0 $ for specificity.

Before delving into the theoretical details, we provide an overview of  eigenfunction estimation in functional data analysis. We start with the resolvent series {shown in Equation \eqref{eq:ptb-hall}} and illustrate its application in statistical analysis,
\begin{equation}\label{eq:ptb-hall}
\E (\|\hat\phi_{j}-\phi_{j} \|^2)\asymp\sum_{k\neq j}\frac{\E[\{\iint(\hat{C}-C)\phi_{j}\phi_{k} \}^2] }{(\lambda_{k}-\lambda_{j})^2}.
\end{equation}
Such expansions can be found in \cite{bosq2000linear}, \cite{dou2012estimation}, and \cite{li2010uniform}; see Chapter 5 in \cite{hsing2015theoretical} for details. Denote by $\eta_j$ the eigengap of $\lambda_j$, that is, $\eta_{j}:=\min_{k\neq j}|\lambda_{k}-\lambda_{j}|$. An basic rough bound for $\E(\|\hat{\phi}_{j}-\phi_{j}\|^2)$ can be derived from Equation \eqref{eq:ptb-hall} and Bessel's inequality, 
\begin{equation}\label{eq:ptb-maintch}
	\E (\|\hat\phi_{j}-\phi_{j} \|^2)\leq\eta_{j}^{-2} \E (\|\hat{C}-C \rhsnorm^2).
\end{equation}
However, this bound is suboptimal for two reasons. First, while $\eta_j$ is bounded away from $0$ for a fixed $j$, the bound implies that the eigenfunctions converge at the same rate as the covariance function. This is counterintuitive since integration usually brings extra smoothness \citep{cai2006prediction}, which typically results in the eigenfunction estimates converging at a faster rate than the two-dimensional kernel smoothing rate of $\|\hat{C}-C\rhsnorm^2$. Second, for $j$ that diverges with the sample size, $\eta_{j}^{-2}\rightarrow\infty$ in the bound cannot be improved. Assuming $\lambda_j\asymp j^{-a}$, the rate obtained by \eqref{eq:ptb-maintch} is slower than $j^{2a+2}/n$, which is known to be  suboptimal \citep{wahl2020information}. Therefore, to obtain a sharp rate for $\E(\|\hat{\phi}_{j}-\phi_{j}\|^2)$, one should use the original perturbation series given by \eqref{eq:ptb-hall}, rather than its approximation given by \eqref{eq:ptb-maintch}.

When each trajectory $X_{i}(t)$ is fully observed for all $t \in [0,1]$, the cross-sectional sample covariance $\hat{C}(s,t) = n^{-1} \sum_{i=1}^{n} X_{i}(s)X_{i}(t)$ is a canonical estimator of $C(s,t)$. Then, the numerators in each term of \eqref{eq:ptb-hall} can be reduced to the principal components scores under some mild assumptions, for example, $\E[\{\iint(\hat{C}-C)\phi_{j}\phi_{k} \}^2] \lesssim n^{-1}\lambda_{j}\lambda_{k}$ \citep{hall2007methodology,dou2012estimation}. Subsequently, $\E( \|\hat\phi_{j}-\phi_{j} \|^2)$ is bounded by $(\lambda_{j}/n)\sum_{k\neq j}\lambda_{k}/(\lambda_{k}-\lambda_{j})^2$. With the common assumption of the polynomial decay of eigenvalues, the aforementioned summation is dominated by $\lambda_{j}/n\sum_{\lfloor j/2 \rfloor\leq k\leq 2j}\lambda_{k}/(\lambda_{k}-\lambda_{j})^2$, which is $O(j^2/n)$ and optimal in the minimax sense \citep{wahl2020information}. See Lemma 7 in \cite{dou2012estimation} for a detailed elaboration. {This suggests that the convergence rate caused by the inverse issue can be captured by the summation over the set $\{k \leq 2j\}$, and the tail sum on $\{k > 2j\}$ can be treated as a unity.}

However, we would like to emphasize that when it comes to discretely observed functional data, all the existing literature utilizing a bound similar to \eqref{eq:ptb-maintch} excludes the case of diverging indices. For instance, the result in \cite{cai2010nonparametric} is simply a direct application of the bound in \eqref{eq:ptb-maintch}. Moreover, their one-dimensional rate is inherited from the covariance estimator, which is assumed to be in a tensor product space that is smaller than the space $\mathcal{L}^2[0,1]^2$. On the other hand, the one-dimensional rates obtained by \cite{hall2006propertieshans} and \cite{li2010uniform} utilize detailed calculations based on the approximation of the perturbation series in \eqref{eq:ptb-hall}. However, these results are based on the assumption that $\eta_{j}$ is bounded away from zero, which implies that $j$ must be a fixed constant. This is inconsistent with the nonparametric nature of functional data  models, which aim to approximate or regularize an infinite-dimensional process. Therefore, when dealing with discretely observed functional data, the key to obtaining a sharp bound for estimated eigenfunctions with diverging indices lies in effectively utilizing the perturbation series  \eqref{eq:ptb-hall}.

The main challenges {arise} from quantifying the summation in \eqref{eq:ptb-hall} without the fully observed sample covariance. For the pre-smoothing method, the reconstructed $\hat{X}_{i}$ achieves a $\sqrt{n}$ convergence in the $\mathcal{L}^{2}$ sense when each $N_{i}$ reaches a magnitude of $n^{5/4}$, and then the estimated covariance function $\hat C(s,t) = n^{-1}\sum_{i=1}^{n}\hat{X}_{i}(s)\hat{X}_{i}(t)$ has an optimal rate $\|\hat C- C\rhsnorm = O_{P}(n^{-1/2})$. However, this does not guarantee optimal convergence of a diverging number of eigenfunctions. The numerators in each term of \eqref{eq:ptb-hall} are no longer the principal component scores, and the complex form of this infinite summation makes it difficult to quantify when $|\lambda_{k}-\lambda_{j}|\rightarrow 0$. Similarly, the pooling method also encounters significant challenges in summing all $\E[\{\iint(\hat{C}-C)\phi_{j}\phi_{k} \}^2]$ with respect to $j$ and $k$. {Specifically, the convergence rate of $\|\hat{C}-C\rhsnorm^2$ should be a two-dimensional kernel smoothing rate with variance $n^{-1}\{1+(Nh)^{-2}\}$ \citep{zhang2016sparse}. However, after being integrated twice, $\E[\{\iint(\hat{C}-C)\phi_{j}\phi_{k} \}^2]$ has a degenerated kernel smoothing rate with variance $n^{-1}$. According to Bessel's inequality, $\E(\|\hat{C}-C\rhsnorm^2)$ can be expressed as $\sum_{j,k}\E[\{\iint(\hat{C}-C)\phi_{j}\phi_{k} \}^2]$. However, due to estimation bias, one cannot directly sum all $\E[\{\iint(\hat{C}-C)\phi_{j}\phi_{k} \}^2]$ with respect to all $j,k$.}

\section{$\mathcal{L}^2$ convergence of eigenfunction estimates}\label{sec:eigthm}
Based on the issues discussed above, we propose a novel technique for addressing the perturbation series \eqref{eq:ptb-hall} when dealing with discretely observed functional data. To this end, we make the following regularity assumptions.

\begin{assum}\label{asm:a1}
	{There exists a positive constant $c_{0}$ such that  $\mathbb{E}(\xi_{j}^{4})\leq c_{0}\lambda_{j}^2$ for all $j$.}
\end{assum}
\begin{assum}\label{asm:cov}
	{The second order derivatives of $C(s,t)$, $\partial C(s,t)/\partial s^2$, $\partial C(s,t)/\partial t^2$ and $\partial C(s,t)/\partial s \partial t$ are bounded on $[0,1]^2$.}
\end{assum}
\begin{assum}\label{asm:a2}
	The eigenvalues $\lambda_{j}$ are decreasing with $ j^{-a}\gtrsim \lambda_{j}\gtrsim\lambda_{j+1} + j^{-a-1}$ for $a >1$ and each $j\geqslant 1$. 
\end{assum}
\begin{assum}\label{asm:a3}
	For each $j\in \mathbb{N}^{+}$, the eigenfunctions $\phi_{j}$ satisfies $\sup\limits_{t\in [0,1]}|\phi_{j}(t)|=O( 1 )$ and
$$\sup\limits_{t\in [0,1]}|\phi_{j}^{(k)}(t)|\lesssim j^{c/2}\sup\limits_{t\in [0,1]}| \phi_{j}^{(k-1)}|, \hspace{0.1in} \text{ for } k=1,2,  $$
where $c$ is a positive constant.
\end{assum}
\begin{assum}\label{asm:u1}
	{$\E [\sup_{s\in[0,1]}|X(t)|^{2\alpha}] <\infty$ and $\E(\varepsilon^{2\alpha})<\infty $ for $\alpha>1$.}
\end{assum}
{Assumptions \ref{asm:a1} and \ref{asm:cov} are widely adopted in the functional data literature related to smoothing \citep{yao2005jasa, cai2010nonparametric,zhang2016sparse}. The decay rate assumption on the eigenvalues provides a natural theoretical framework for justifying and analyzing the properties of functional principal components \citep{hall2007methodology, dou2012estimation, zhou2023functional}.} The number of eigenfunctions that can be well estimated from exponentially decaying eigenvalues is limited to the order of $\log n$, which lacks practical interest. Consequently, we adopt the assumption of polynomial decay in eigenvalues.  To achieve quality estimates for a specific eigenfunction, its smoothness should be taken into account. Generally, the frequency of $\phi_{j}$ is higher for larger $j$, which requires a smaller bandwidth to capture its local variation. Assumption \ref{asm:a3} characterizes the frequency increment of a specific eigenfunction via the smoothness of its derivatives. For some commonly used bases, such as the Fourier, Legendre, and wavelet bases, \( c=2 \). In \cite{hall2006propertieshans}, the authors assumed that \(\max_{1\leq j\leq r}\max_{s=0,1,2}\sup_{t\in[0,1]}|\phi_{j}^{(s)}(t) |\leq C\), which is only achievable for a fixed \( r \). Therefore, Assumption \ref{asm:a3} can be interpreted as a generalization of this assumption. {To analyze the convergence of eigenfunctions effectively, a uniform convergence rate of the covariance function is needed to handle the local linear estimator (\ref{eq:ll}). Assumption \ref{asm:u1} is the moment assumption required for uniform convergence of covariance function and adopted in \cite{li2010uniform} and \cite{zhang2016sparse}.}

{For the observation time points $\{t_{ij}\}_{i,j}$, there are two typical types of designs: the random design, in which the design points are random samples from a distribution, and the regular design, where the observation points are predetermined  mesh grid. For the random design, the following assumption is commonly adopted \citep{yao2005jasa, li2010uniform, cai2011optimal, zhang2016sparse}:}
\begin{assum}[Random design]\label{asm:r}
	{The design points $t_{ij}$, which are independent of $X$ and $\varepsilon$, are i.i.d. sampled from a distribution on $[0,1]$ with a density that is bounded away from zero and infinity.}
\end{assum}
{For the regular design, each sample path is observed on an equally spaced grid $\{t_{j}\}_{j=1}^{N}$, where $N_{i}=\cdots = N_{n} = N$ for all subjects. This longitudinal design is frequently encountered in a various scientific experiments and has been studied in \cite{cai2011optimal}. Assumption \ref{asm:e} guarantees a sufficient number of observations within the local window for the kernel smoothing method. Furthermore, Assumption \ref{asm:e-c} is needed for the  Riemann sum approximation in the fixed regular design.}
\begin{assum}[Fixed regular design]\label{asm:e}
	{The design points $\{t_{j}\}_{j=1}^{N}$ are nonrandom, and there exist constant $c_{2}\geq c_{1}>0 $, such that for any interval $A,B\in [0,1]$,
	\begin{itemize}
		\item[(a)] $c_{1}N|A|-1\leq \sum_{j=1}^{N}\mathds{1}_{\{t_{j}\in A\}}\leq \max\{c_{2}N,1 \} $,
		\item[(b)] $c_{1}N^2|A||B|-1\leq\sum_{l_1,l_2}^{N} \mathds{1}_{\{t_{l_1}\in A\}} \mathds{1}_{\{t_{l_2}\in B\}}\leq \max\{c_2 N^2 |A||B|,1 \}$, 
	\end{itemize}
	where $|A|$ denotes the length of $A$.} 
\end{assum}
\begin{assum}\label{asm:e-c}
	${\partial C(s,t)}/{\partial s^2}$ and ${\partial C(s,t)}/{\partial t^2}$ are continuously differentiable.
\end{assum}
{The following theorem is one of our main results. It characterizes the $\mathcal{L}^{2}$ convergence  of the estimated eigenfunctions with diverging indices for both random design (Assumption \ref{asm:r}) and fixed regular design (Assumption \ref{asm:e}).}
{\begin{theorem}\label{thm-eig}
	Assume Assumptions \ref{asm:a1} to \ref{asm:a3} hold, further assume Assumption \ref{asm:u1} holds with $\alpha>3$.
	\begin{itemize}
		\item[(a)] For the random design, under Assumption \ref{asm:r}, for all $j\leq m\in\mathbb{N}_{+}$ satisfies $m^{2a+2}/n\rightarrow0$, $m^{2a+2}/(n\bar{N}_2^2h^2)\rightarrow0 $ and $h^{4}\max\{ m^{2a+2c},m^{4a}\log n \}\lesssim1 $,
		\begin{equation}\label{eq:eig}
		\|\hat{\phi}_{j}-\phi_{j}\|^2=O_{P} \left(\frac{j^{2}}{n}\left\{1 +\frac{j^{2a}}{\bar{N}_2^2}  \right\}+\frac{j^{a}}{n\bar{N}_2h}\left(1 + \frac{j^{a}}{\bar{N}_2}\right) +h^{4}j^{2c+2}\right),
		\end{equation}
		where $\bar{N}_{2}=(n^{-1}\sum_{i=1}^{n}N_{i}^{-2})^{-1/2}$.
		\item[(b)] For the fixed regular design, under Assumption \ref{asm:e} and \ref{asm:e-c}, for all $j\leq m\in\mathbb{N}_{+} $ satisfies $m^{2a+2}/n\rightarrow0$, $m^{2a+2}/(n{N}^2h^2)\rightarrow0$, ${N}h^{a}\gtrsim 1$ and $h^{4}\max\{ m^{2a+2c},m^{4a}\log n \}\lesssim1 $,
		\begin{equation}\label{eq:eig-f}
			\|\hat{\phi}_{j}-\phi_{j}\|^2=O_{P} \left(\frac{j^{2}}{n}+\frac{j^a}{nNh} +h^{4}j^{2c+2}\right).
		\end{equation}
	\end{itemize}
\end{theorem}}

The integer \( m \) in Theorem \ref{thm-eig} represents the maximum number of eigenfunctions that can be well estimated from the observed data using appropriate tuning parameters. It is important to note that in Theorem \ref{thm-eig}, \( m \) could diverge to infinity. The upper bound of \( m \) is a function of the sample size $n$, the sampling rate $\bar{N}_2$, the smoothing parameter $h$, and the decaying eigengap $\eta_{j}$ or $a$. As the frequency of $\phi_{j}$ is higher for larger $j$, smaller $h$ is required to capture its local variations. If \( a \) is large, the eigengap \( \eta_{j} \) diminishes rapidly, posing a greater challenge in distinguishing between adjacent eigenvalues. {Note that the assumptions of $m$ in Theorem \ref{thm-eig} could encompass most downstream analyses that involve a functional covariate, such as functional linear regression as discussed in \citep{ hall2007methodology}.}

{Theorem \ref{thm-eig} provides a good illustration of both the infinite dimensionality and smoothness nature of functional data. To better understand this result, note that \( j^2/n \) is the optimal rate in the fully observed case. The additional terms on the right-hand side of Equation \eqref{eq:eig} represent contamination introduced by discrete observations and measurement errors. In particular, the term with \( h^4 \) corresponds to the smoothing bias, while the term \( (n\bar{N}_2h)^{-1} \) reflects the variance typically associated with one-dimensional kernel smoothing. Terms including \( \bar{N}_2^{-1} \) arise from the discrete approximation, and the terms that involve \( j \) with positive powers are due to the decaying eigengaps associated with an increasing number of eigencomponents.}

{The first part of Theorem \ref{thm-eig} provides a unified result for eigenfunction estimates under random design without imposing any restrictions on the sampling rate \(N_i\). Similar to the phase transitions of mean and covariance functions studied in \cite{cai2011optimal} and \cite{zhang2016sparse}, Corollary \ref{cor:eig} presents a systematic partition that ranges from ``sparse'' to ``dense'' schemes for eigenfunction estimation under the random design, which is meaningful for FPCA-based models and methods.}

\begin{corollary}\label{cor:eig}
Under Assumptions \ref{asm:a1} to \ref{asm:r}, and $m\in\mathbb{N}_{+}$ satisfies $(a)$ of Theorem \ref{thm-eig}. For each $j\leq m$ and let $h_{opt}(j)=(n\bar{N}_{2})^{-1/5}j^{(a-2c-2)/5}(1+j^{a}/\bar{N}_2)^{1/5} $,
	\begin{enumerate}[label=\textup{(\alph*)}]
		\item If $\bar{N}_2\gtrsim  j^{a}$,
		$$
			\|\hat\phi_{j}-\phi_{j}\|^2=O_{P}\left(  \frac{j^{2}}{n}+\frac{j^{(4a+2c+2)/5}}{(n\bar{N}_2)^{4/5}}\right) .
		$$ 
		In addition, if $\bar{N}_2\geq n^{1/4}j^{a+c/2-2}$, $$ \|\hat\phi_{j}-\phi_{j}\|^2=O_{P}\left(\frac{j^2}{n} \right) .$$
		\item If $\bar{N}_2=o(j^{a})$, 
		$$
			\|\hat\phi_{j}-\phi_{j}\|^2=O_{P}\left(  \frac{j^{2a+2}}{n\bar{N}_2^2}+\frac{j^{(8a+2c+2)/5}}{(n\bar{N}_2^{2})^{4/5}}\right). 
		$$
	\end{enumerate}
\end{corollary}

Note that $h_{opt}(j)$ is the optimal bandwidth for estimating a specific eigenfunction $\phi_{j}$. However, in practice, it suffices to estimate the covariance function just once using the optimal bandwidth associated with the largest eigenfunction \( \phi_{m_{\text{max}}} \). This  ensures that both the subspace spanned by the first \( m_{\text{max}} \) eigenfunctions and their corresponding projections are well estimated. {When conducting downstream analysis with FPCA, the evaluation of error rates typically involves the term $\sum_{j=1}^{m_{\max}} j^{-\beta}\|\hat{\phi}_j-\phi_j\|^2$. Here $m_{\max}$ denotes the maximum number of principal components  used, and $\beta$ varies across different scenarios, reflecting how the influence of each principal component on the error rate is weighted. For example, in the functional linear model,   \(\beta>0\) represents the rate at which the Fourier coefficients of the regression function decay. By the first statement of Theorem \ref{thm-eig}, 
$$
		\begin{aligned}
			\sum_{j=1}^{m_{\max}}j^{-\beta}\|\hat\phi_j-\phi_j\|^{2}=\frac{m_{\max}^{3-\beta}}{n}\left(1+\frac{m_{\max}^{2a}}{\bar{N}_2^2} \right)+\frac{m^{a-\beta+1}}{n\bar{N}_2h}\left(1+\frac{m_{\max}^{a}}{\bar{N}_2} \right)+h^{4}m^{2c+3-\beta}.
		\end{aligned}
$$
It can be found that the optimal bandwidth $h$ in the above equation aligns with $h_{opt}(m_{\max})$ for all $\beta$. This implies that our theoretical framework can be seamlessly adapted to scenarios where a single bandwidth is employed to attain the optimal convergence rate in the downstream analysis.}

{For the commonly used bases where $c=2$, the convergence rate for the $j$th eigenfunction achieves optimality as if the curves are fully observed when $\bar{N}_2> n^{1/4}j^{a-1}$. Keeping $j$ fixed, the phase transition occurs at $n^{1/4}$, aligning with results  in \cite{hall2006propertieshans} and \cite{cai2010nonparametric}, as well as mean and covariance functions discussed in \cite{zhang2016sparse}. For $n$ subjects, the maximum index of the eigenfunction that can be well estimated is  less than  $m_{\mathrm{max}}:=n^{1/(2a+2)}$, which directly follows from the assumption $m^{2a+2}/n\rightarrow 0$.The phase transition in estimating $\phi_{m_{\mathrm{max}}}$ occurs at $n^{1/4 + (a-1)/(2a+2)}$.} This can be interpreted from two aspects. On one hand, compared to mean and covariance estimation, more observations per subject are required to obtain optimal convergence for eigenfunctions with increasing indices, showing the challenges tied to infinite dimensionality and decaying eigenvalues. On the other hand, the fact that $n^{1/4 + (a-1)/(2a+2)}$ is only slightly larger than $n^{1/4}$ justifies the merits of the pooling method and supports our intuition. When $j$ is fixed and $\bar{N}_2$ is finite, the convergence rate obtained by Corollary \ref{cor:eig} is $(nh)^{-1} + h^4$, which corresponds to a typical one-dimensional kernel smoothing rate and achieves optimal at $h \asymp n^{-1/5}$. This result aligns with those in \cite{hall2006propertieshans} and is optimal in the minimax sense. When allowing $j$ to diverge, the known lower bound  $j^2/n$ for fully observed data is attained by applying van Trees' inequality to the special orthogonal group \citep{wahl2020information}. However, the argument presented in \cite{wahl2020information} does not directly extend to the discrete case and there are currently no available lower bounds for the eigenfunctions with diverging indices based on discrete observations, which remains an open problem that requires further investigation.

%Note that $h_{opt}(j)$ is the optimal bandwidth for estimating a specific eigenfunction $\phi_{j}$. In the theoretical analysis, higher-order eigenfunctions \( \phi_{j} \) exhibit increased frequency for larger values of \( j \), necessitating a smaller \( h \) to accurately capture their local variations. Thus, the optimal bandwidth should  dependent on the index $j$ theoretically in the context of analyzing the convergence rate of a specific eigenfunction. \hang{In practice, it suffices to estimate the covariance function just once using the optimal bandwidth associated with the largest eigenfunction, \( \phi_{m_{\text{max}}} \). This approach ensures that both the subspace spanned by the first \( m_{\text{max}} \) eigenfunctions and their corresponding projections are well estimated. Take functional linear regression as an example: here, the regression estimator is constructed using the first \( m_{\text{max}} \) eigenfunctions and eigenvalues. The theoretically optimal bandwidth for the regression function is \( h_{\text{opt}}(m_{\text{max}}) \). In downstream analysis, encountering a subspace generated by the first \( m \) eigenfunctions is more typical than the application of a specific eigenfunction. Our method requires only a single bandwidth to optimally generate the subspace formed by the first $m$ eigenfunctions.}

Comparing the results of this work with those in \cite{zhou2023functional} is also of interest. The $\mathcal{L}^2$  convergence rate for the $j$th eigenfunction obtained  in \citep{zhou2023functional} is
	\begin{equation}\label{eq:flr}
		\frac{j^{a+2}}{n}\left(1+\frac{1}{Nh} \right)+h^{4}j^{2a+2c+2}.
	\end{equation} It is evident that the results in \cite{zhou2023functional} will never reach the optimal rate  $j^{2}/n$ for any sampling rate $N$. In contrast, Corollary \ref{cor:eig}  provide a systematic partition ranging from ``sparse'' to ``dense'' schemes for eigenfunction estimation, and the optimal rate  $j^{2}/n$ can be achieved when the sampling rate $N_i$ exceeds the phase transition point. {The optimal rate achieved here represents more than just a theoretical improvement over previous findings; it also carries substantial implications for downstream analysis. Further discussion can be found in Section \ref{sec:conc}.}

%TODO fix design and unknown mean
{The following Corollary presents the phase transition of eigenfunctions for fixed regular design. Note that for fixed regular design, the number of distinct observed time points in an interval of length $h$ is on the order of $Nh$, so $Nh \gtrsim 1$ is required to ensure there is at least one observation in the bandwidths of kernel smoothing \citep{shao2022intrinsic}. Moreover, $Nh^a\gtrsim 1 $ is required to eliminate the Riemann sum approximation bias. The condition $Nh^a\gtrsim 1 $ is parallel of part (a) in Corollary \ref{cor:eig} in the scenario of the random design, which is similar as the mean and covariance estimation where consistency can only be obtained under the dense case for regular design \citep{shao2022intrinsic}.}

{\begin{corollary}\label{cor:eig-fix}
Under Assumptions \ref{asm:a1} to \ref{asm:u1}, \ref{asm:e}, \ref{asm:e-c} and $m\in\mathbb{N}_{+}$ satisfies $(b)$ of Theorem \ref{thm-eig}. For each $j\leq m$ and let $h_{opt}(j)= (nN)^{-1/5}j^{(a-2c-2)/5} $,
$$
			 \|\hat\phi_{j}-\phi_{j}\|^2=O_{P}\left( \frac{j^{2}}{n}+\frac{j^{(4a+2c+2)/5}}{(n{N})^{4/5}}\right) .
		$$
		In addition, if ${N}\geq n^{1/4}j^{a+c/2-2}$, $$\|\hat\phi_{j}-\phi_{j}\|^2 =O_{P} \left( \frac{j^2}{n} \right) .$$
\end{corollary}}

If the mean function $\mu(t)$ is unknown, one could use local linear smoother to fit $\hat{\mu}(t)=\hat{\alpha}_{0}$ with
$$(\hat{\alpha}_{0},\hat{\alpha}_{1})=\underset{\alpha_0, \alpha_1}{\operatorname{argmin}}\frac{1}{n} \sum_{i=1}^{n} \frac{1}{N_i} \sum_{j=1}^{N_i}\{X_{ij}-\alpha_{0}-\alpha_{1}(t_{ij}-t) \}^2\frac{1}{h_\mu} \K\left(\frac{t_{ij}-t}{h_{\mu}} \right) .$$
Then the covariance estimator $\hat{C}$ is obtained by replacing the raw covariance $\delta_{ijl}$ in \eqref{eq:ll} by $\hat{\delta}_{ijl}=\{X_{ij}-\hat{\mu}(t_{ij})\}\{X_{il}-\hat{\mu}(t_{il})\}$.  The following corollary presents the convergence rate and  phase transition for the case where $\mu(t)$ is unknown.
It should be noted that the Fourier coefficients of $\mu(t)$ with respect to the eigenfunction ${\phi_{j}}$ generally do not exhibit a decaying trend. Therefore, to eliminate the estimation error caused by the mean estimation, an additional lower bound on $N$ is necessary. This lower bound, denoted as $N \gtrsim j^{a}$, aligns with the partition described in Corollary \ref{cor:eig} for the random design case.  
\begin{corollary}\label{cor:eig_umean}
	Suppose that Assumptions \ref{asm:a1} to \ref{asm:a3} hold.  Under either of Assumption  \ref{asm:r} or \ref{asm:e} and \ref{asm:e-c}, for all  $m\in\mathbb{N}_{+}$ satisfies $m^{2a+2}/n$, $Nh^{a}\gtrsim1$ and $h^{4}\max\{m^{2a+2c},m^{4a}\log n \}\lesssim1$ and each $j\leq m$,
$$
			 \|\hat\phi_{j}-\phi_{j}\|^2=O_{P}\left( \frac{j^{2}}{n}+\frac{j^{(4a+2c+2)/5}}{(n\bar{N}_2)^{4/5}}\right) .
		$$
		In addition, if $\bar{N}_2\geq n^{1/4}j^{a+c/2-2}$, $$\|\hat\phi_{j}-\phi_{j}\|^2 =O_{P} \left( \frac{j^2}{n} \right) .$$
\end{corollary}

\section{Asymptotic normality of eigenvalue estimates}\label{sec:asy}

The distribution of eigenvalues plays a crucial role in statistical learning and is of significant interest in the high-dimensional setting. Random matrix theory provides a systematic tool for deriving the distribution of eigenvalues of a squared matrix \citep{anderson2010introduction, pastur2011eigenvalue}, and has been successfully applied in  various statistical problems, such as signal detection \citep{nadler2011performance, onatski2009testing, bianchi2011performance}, spiked covariance models \citep{johnstone2001distribution, paul2007asymptotics, el2007tracy, ding2021spiked, bao2022statistical}, and hypothesis testing \citep{bai2009corrections, chen2010two, zheng2012central}. For a comprehensive treatment of random matrix theory in statistics, we recommend the monograph by \cite{bai2010spectral} and the review paper by \cite{paul2014random}.

Despite the success of random matrix theory in high-dimensional statistics, its application to functional data analysis is not straightforward due to the different  structures of functional spaces. If the observations are taken at the same {time points} $\{t_{j}\}_{j=1}^{N}$ for all $i$,  one can obtain an estimator for $\Sigma_T =\cov (\tilde{X}_i,\tilde{X}_i)+\sigma_{X}^{2}I_{N} $, where $\tilde{X}_i=(X_{i}(t_1),\ldots,X_{i}(t_N) )^{\top} $. Note that $\tilde{X}_i$ can be regarded as a random vector; {however} the adjacent elements in $\tilde{X}_i$ are highly correlated as $N$ increases due to the smooth nature of functional data. {This correlation violates the independence assumption required in most random matrix theory settings.}

In the context of functional data, variables of interest become the eigenvalues of the covariance operator. However, the rough bound obtained by Weyl's inequality, $|\hat\lambda_{k}-\lambda_{k}|\leq \|\Delta\rhsnorm$, is suboptimal {from two respects}. First, $|\hat\lambda_{k}-\lambda_{k}|$ should have a degenerated kernel smoothing rate with variance $n^{-1/2}$, whereas $\|\Delta\rhsnorm$ has a two-dimensional kernel smoothing rate with variance $(n\bar{N}_2^2h^2)^{-1/2}$. Second, due to the infinite dimensionality of functional data, the eigenvalues $\{\lambda_{k}\}_k$ tend to zero as $k\rightarrow\infty$, so a general bound for all eigenvalues provides little information for those with larger indices. {Although expansions and asymptotic normality have been studied for a fixed number of eigenvalues, as well as for those with diverging indices for fully observed functional data, the study of eigenvalues with a diverging index under the discrete sampling scheme remains limited.}

{In light of the aforementioned issues, we employ the perturbation technique outlined in previous sections to establish the asymptotic normality of eigenvalues with diverging indices, which holds broad implications for inference problems in functional data analysis. Before presenting our results, we introduce the following assumption, which is standard in FPCA for establishing asymptotic normality.}

\begin{assum}\label{asm:n1}
	$\E(\|X\|^{6})<\infty$ and $\E(\epsilon^{6})<\infty$. For any sequence $j_{1},\ldots,j_{4}\in\mathbb{N}_{+}$, $\E (\xi_{j_{1}}\xi_{j_{2}}\xi_{j_{3}}\xi_{j_{4}} )=0$ unless each index $j_{k}$ is repeated.
\end{assum}

%
%Assumption \ref{asm:n1} introduces the  moment conditions necessary to establish the asymptotic normality, as detailed in \citep{zhang2016sparse}. This assumption is a standard prerequisite in the context of Functional Principal Component Analysis (FPCA), referenced in \citep{cai2006prediction, hall2009theory}, and serves to streamline the computational aspects.

\begin{theorem}\label{thm:egv}
Under Assumptions \ref{asm:a1} to \ref{asm:a3}, \ref{asm:r} and \ref{asm:n1},  for all $j\leq m\in\mathbb{N}_{+}$ satisfy $m=o(n^{1/(2a+4)})$, $h{m^{a}}\log n \lesssim1$ and $h^{4}m^{2a+2c}\lesssim1$
	$$\Sigma_{n}^{-\frac{1}{2}}\left(\frac{\hat\lambda_{j}-\lambda_{j}}{\lambda_{j}}-{\K_2 h^{2}}\int\phi_{j}^{(2)}(u)\phi_{j}(u)\diff u+o(h^{2}) \right)\stackrel{d}{\longrightarrow}\mathcal{N}(0,1), $$
		where
		\begin{align*}
	&\Sigma_{n}=\frac{4!P_{0}(N)}{n}\frac{\E(\xi_{j}^4)-\lambda_{j}^2}{\lambda_{j}^2}+4!\frac{P_{1}(N)}{n} \frac{\int \{C(u,u)+\sigma_{X}^2 \}\frac{\phi_{j}^2(u)}{f(u)}\diff u }{\lambda_{j}}\\
	+&4\frac{P_2(N)}{n}\frac{1}{\lambda_{j}^2} \left(\left[ \int\{C(u,u)+\sigma_{X}^2 \}\frac{\phi_{j}^2(u)}{f(u)}\diff u \right]^2-\iint C(u,v)^2\frac{\phi_{j}^2(u)\phi_{j}^2(v)}{f(u)f(v)}\diff u \diff v \right)
\end{align*}
with $\K_2=\int u^2K(u)\diff u$ and
$$P_{0}(N)=\frac{1}{n}\sum_{i=1}^{n}\frac{(N_i-2)(N_i-3)}{N_i(N_i-1)},\ P_{1}(N)=\frac{1}{n}\sum_{i=1}^{n}\frac{(N_i-2)}{N_i(N_i-1)}, \ P_{2}(N)=\frac{1}{n}\sum_{i=1}^{n}\frac{1}{N_i(N_i-1)}.$$
\end{theorem}

Since $\lambda_j$ approaches zero as $j$ approaches infinity, we need to regularize the eigenvalues so that they can be compared on the same scale of variability. {For a fixed $j$, $(\hat{\lambda}_j - \lambda_j)/\lambda_j$ is $\sqrt{n}$-consistent when $h$ is small. For diverging $j$, Corollary \ref{cor-egv} presents three different types of asymptotic normalities, depending on the value of $P_{2}(N)$.} 

\begin{corollary}\label{cor-egv}
	Under the assumptions of Theorem \ref{thm:egv}, 
	\begin{enumerate}[label=\textup{(\alph*)}]
	\item If $P_{2}(N)/\lambda_{j}^2\rightarrow0 $, $\sqrt{n}h^{2}\int\phi_{j}^{(2)}(u)\phi_{j}(u)\diff u\rightarrow 0 $,
	$$\sqrt{n}\left( \frac{\hat\lambda_{j}-\lambda_{j}}{\lambda_{j}} \right)\stackrel{d}{\longrightarrow}\mathcal{N}\left(0,4! \frac{\E (\xi_{j}^{4})-\lambda_{j}^{2}}{\lambda_{j}^{2}} \right).  $$
	\item If $P_{2}(N)/\lambda_{j}^2\rightarrow C_3 $ for a positive $C_3$,
	\begin{align*}
		&\sqrt{n}\left( \frac{\hat\lambda_{j}-\lambda_{j}}{\lambda_{j}}-{\K_2 h^{2}}\int\phi_{j}^{(2)}(u)\phi_{j}(u)\diff u \right)
		\stackrel{d}{\longrightarrow}\mathcal{N}\left(0, 4!\frac{\E( \xi_{j}^{4})-\lambda_{j}^{2}}{\lambda_{j}^{2}}\right.\\
		&\quad\quad +4!\sqrt{C_3} \int \{C(u,u)+\sigma_{X}^2 \}\frac{\phi_{j}^2(u)}{f(u)}\diff u \\
		&\quad\quad+\left.4C_3 \left\{\left[ \int\{C(u,u)+\sigma_{X}^2 \}\frac{\phi_{j}^2(u)}{f(u)}\diff u \right]^2-\iint C(u,v)^2\frac{\phi_{j}^2(u)\phi_{j}^2(v)}{f(u)f(v)}\diff u \diff v \right\}\right)
	\end{align*}
	\item If $P_{2}(N)/\lambda_{j}^2\rightarrow\infty $,
	\begin{align*}
		&\sqrt{\frac{n\lambda_{j}^2}{P_{2}(N)}}\left( \frac{\hat\lambda_{j}-\lambda_{j}}{\lambda_{j}}-{\K_2 h^{2}}\int\phi_{j}^{(2)}(u)\phi_{j}(u)\diff u\right)\\
		&\stackrel{d}{\longrightarrow}\mathcal{N}\left(0, 4\left\{\left[ \int\{C(u,u)+\sigma_{X}^2 \}\frac{\phi_{j}^2(u)}{f(u)}\diff u \right]^2-\iint C(u,v)^2\frac{\phi_{j}^2(u)\phi_{j}^2(v)}{f(u)f(v)}\diff u \diff v \right\} \right).
	\end{align*}
	\end{enumerate}

\end{corollary} 

Compared to the mean and covariance estimators, which are associated with one-dimensional and two-dimensional kernel smoothing rates respectively, the estimator of eigenvalues exhibits a degenerate rate after being integrated twice. This implies that there is no trade-off between bias and variance in the bandwidth $h$, and the estimation bias can be considered negligible for small values of $h$. In this scenario, the phase transition is entirely determined by the relationship between $P_{2}(N)$ and $\lambda_{j}$. Specifically, in the dense and ultra-dense cases where $P_{2}(N)/\lambda_{j}^2<\infty $, the variance terms  resulting from discrete approximation are dominated by $1/\sqrt{n}$, corresponding to cases (a) and (b) in Corollary \ref{cor-egv}. On the other hand, when each $N_i$ is relatively small, the estimation variance $\sqrt{n\lambda_{j}^2/P_{2}(N)}$ arising from discrete observations dominates, as outlined in case (c) of Corollary \ref{cor-egv}. The phase transition point $j^{a}$ is the same as the case of eigenfunctions outlined in Corollary \ref{cor:eig} and \ref{cor-eigunif-diverg}, and for larger values of $j$, more observations are needed due to the vanishing eigengap.

{Using similar arguments, we establish the  $\mathcal{L}^2$ convergence rate for eigenvalues, which is useful for analyzing models involving functional principal component analysis, such as the plug-in estimator in \cite{hall2007methodology}.}

\begin{corollary}\label{cor-ev}
	Under Assumptions \ref{asm:a1} to \ref{asm:a3} and \ref{asm:r},  for all $j\leq m\in\mathbb{N}_{+}$ satisfy $m=o(n^{1/(2a+4)})$, $h{m^{a}}\log n \lesssim1$ and $h^{4}m^{2a+2c}\lesssim1$,
$$(\hat{\lambda}_j-\lambda_j )^{2}=O_{P}\left(\frac{j^{-2a}}{n}\left\{P_{0}(N)+j^{2a}P_{2}(N) \right\}+h^4 j^{2c-2a}   \right) .$$
\end{corollary}

\section{Uniform convergence of covariance and eigenfunction estimates}\label{sec:uni}

Classical nonparametric regression with independent observations has yielded numerous results for {the} uniform convergence of kernel-type estimators \citep{bickel1973some, hardle1988strong, hardle1989asymptotic}. For functional data with in-curve dependence, \cite{yao2005jasa} obtained a  uniform bound for mean and covariance pooling estimates. More recently, \cite{li2010uniform} and \cite{zhang2016sparse} {have} studied {the} strong uniform convergence {of} these estimators, showing that these rates depend on both the sample size and the number of observations per subject. However, {uniform results  for estimated eigenfunctions with diverging indices have not been obtained, even in the fully observed case.}

{Even for the covariance estimates, there remains a theoretical challenge in achieving  uniform convergence for the covariance function under the dense/ultra-dense schemes. Specifically,} to obtain uniform convergence with in-curve dependence in functional data analysis, a common {approach} involves employing the Bernstein inequality to obtain a uniform bound over a finite grid of the observation domain. {This grid becomes increasingly dense as the sample size grows. The goal is then to demonstrate that the bound over the finite grid and the bound over $[0,1]$ are asymptotically equivalent.} This technique has been studied by \cite{li2010uniform} and \cite{zhang2016sparse} to {achieve} uniform convergence for mean and covariance estimators based on local linear smoothers. However, due to the {lack of compactness in functional data}, truncation on the observed data  is necessary to apply the Bernstein inequality; that is,  $X_{ij}\mathds{1}_{\{X_{ij}\leq A_{n}\}}$, where $A_{n}\rightarrow\infty$ as $n\rightarrow\infty$. The choice of the truncation sequence $A_n$ should balance the trade-off between the estimation variance appearing in the Bernstein inequality and the bias resulting from truncation. Once  the optimal $A_n$ is chosen, {it is essential to impose additional moment conditions on both $X_i$ and $\varepsilon_i$ to ensure that the truncation bias is negligible.} For covariance estimation, {the current state-of-the-art results \citep{zhang2016sparse}} require assuming that the $6/(1-\tau)$th moments of $X_i$ and $\varepsilon_i$ are finite, where $\tau$ is the sampling rate with $\bar{N}_{2}\asymp n^{\tau/4}$ for $\tau\in[0,1)$. However, when $\bar{N}_{2}$ is close to $n^{1/4}$, the value of $6/(1-\tau)$ {tends} to infinity. Moreover, in the dense and ultra-dense sampling schemes where $\bar{N}_2\gtrsim n^{1/4}$, {the discrepancy between the truncated and the original estimators becomes dominant, preventing the achievement of optimal convergence rates with the current methods.}

{We first resolve this issue for covariance estimates. We propose an additional truncation on the summation of random quantities after truncating $X_{ij}$ on $A_{n}$, which achieves a sharp bound by Bernstein inequality twice and allows for a larger $A_n$ to reduce the first truncation bias. This novel double truncation technique enables  one to obtain   a unified result for all sampling schemes and to address the limitations of the existing literature.}
% To proceed,  we shall make the following assumption,  which is an additional moment condition required for uniform convergence and adopted  by \cite{li2010uniform, zhang2016sparse}. }

\begin{theorem}\label{thm-Cunif}
	Under the Assumptions \ref{asm:a1}, \ref{asm:u1} and \ref{asm:r},
	\begin{enumerate}[label=\textup{(\alph*)}]
		\item For $p,q=0,1$, denote $$R_{pq}(s,t)=\sum_{i=1}^{n}v_{i}\sum_{l_{1}\neq l_{2}}^{N_{i}}\frac{1}{h^2}\K\left(\frac{t_{il_1}-s}{h} \right) \K\left(\frac{t_{il_2}-t}{h} \right)\left(\frac{t_{il_1}-s}{h} \right)^p \left(\frac{t_{il_2}-t}{h} \right )^{q}\delta_{il_1l_2}.
 $$ Then \begin{equation}\label{eq:Cunif}
		\begin{aligned}
			&\sup_{s,t\in[0,1]}|{R}_{pq}(s,t)- \E{R}_{pq}(s,t)|\\
			=&O\left(  \sqrt{\frac{\ln n}{n}}\left(1+\frac{1}{\bar{N}_2h} \right)+
\left|\frac{\ln n}{n}\right|^{1-\frac{1}{\alpha}}\left|1+\frac{\ln n}{\bar{N}_{2}h}\right|^{2-\frac{2}{\alpha}}h^{-\frac{2}{\alpha}}\right) \text{ a.s. }
		\end{aligned}
	\end{equation}
		\item In addition, if Assumption \ref{asm:cov} holds and   $\alpha>3$,
		\begin{itemize}
		\item[(1).] When $\bar{N}_2/(n/\log n)^{1/4}\rightarrow0 $ and  $h\asymp(n\bar{N}_2^2/\log n)^{-1/6} $,
		$$ \sup_{s,t\in[0,1]}|\hat{C}(s,t)- {C}(s,t)| \lesssim  \left(\frac{\log n}{n\bar{N}_2^2} \right)^{\frac{1}{3}} \quad {a.s.}  $$ 
		\item[(2).] When $\bar{N}_2/(n/\log n)^{1/4}\rightarrow C_1 $ for a positive constant $C_1$ and $h\asymp(n/\log n)^{1/4} $
		$$ \sup_{s,t\in[0,1]}|\hat{C}(s,t)- {C}(s,t)| \lesssim  \sqrt{\frac{\log n}{n}} \quad {a.s.}  $$
		\item[(3).] When $\bar{N}_2/(n/\log n)^{1/4}\rightarrow\infty$,  $h=o(n/\log n)^{1/4} $ and $\bar{N}_2h\rightarrow\infty$,
		$$ \sup_{s,t\in[0,1]}|\hat{C}(s,t)- {C}(s,t)| \lesssim  \sqrt{\frac{\log n}{n}} \quad {a.s.}  $$
	\end{itemize}
	\end{enumerate}
		
\end{theorem}

The first statement of Theorem \ref{thm-Cunif} establishes the uniform convergence rate of the variance term of $\hat{C}-C $ for all sampling rates $\bar{N}_2$. {It is worth comparing} our results with those {in} \cite{li2010uniform} and \cite{zhang2016sparse}. The bias term caused by double truncation, which appears in the second term on the right-hand side of \eqref{eq:Cunif}, is smaller than those obtained in \cite{li2010uniform} and \cite{zhang2016sparse}. As a result, the second statement of Theorem \ref{thm-Cunif} shows that the truncation bias is dominated by the main term if $\alpha>3$, a condition that is much milder compared to the moment assumptions in \cite{li2010uniform} and \cite{zhang2016sparse}, {thereby  establishing} the uniform convergence rate for both sparse and dense functional data. Moreover, when $1<\alpha\leq 3$, {which corresponds to the case where  $X(t)$ or $\varepsilon$ do not possess a sixth order finite moment, the additional term caused by truncation becomes dominant.} In summary, by introducing the double truncation technique, we resolve {the aforementioned issues present in the original proofs of \cite{li2010uniform} and \cite{zhang2016sparse}}, achieve the uniform convergence rate for the covariance estimator {across} all sampling schemes, and show that the optimal rates for dense functional data can be obtained as a special case. Having set the stage, we arrive at the following result that gives the uniform convergence for estimated eigenfunctions with diverging indices.

\begin{theorem}\label{thm-eigunif}
	Under Assumptions \ref{asm:a1} to \ref{asm:u1}, for all $j\leq m$ satisfying $m^{2a+2}/n=o(1)$, $m^{2a+2}/(n\bar{N}_2^2h^2)=o(1)$, $h^{4}m^{2a+2c}\lesssim1$ and $hm^a\log n\lesssim1$,
	\begin{equation}\label{eq:eigunif}
	\begin{aligned}
		\|\hat\phi_{j}-\phi_{j}\|_{\infty}=&O\left( \frac{j}{\sqrt{n}}(\sqrt{\ln n}+\ln j){\left\{ 1+\frac{j^a}{\bar{N}_2}+\sqrt{\frac{j^{a-1}}{\bar{N}_2h}}\left(1+\sqrt{\frac{j^{a}}{\bar{N}_2}} \right) \right\}}\right.\\
			&+\left.j^a\left|\frac{\ln n}{n}\right|^{1-\frac{1}{\alpha}}\left|j^{1/2}+\frac{\ln n}{\bar{N}_2h}\right|^{1-\frac{1}{\alpha}}h^{-\frac{1}{\alpha}}+h^{2}j^{c+1}\log j\right).	\end{aligned}
	\end{equation}
\end{theorem}

The contributions of Theorem \ref{thm-eigunif} are two-fold. First, our result is the first to establish uniform convergence for eigenfunctions with diverging indices,   providing a useful tool for the theoretical analysis of models involving FPCA and inverse issues. Second, the double truncation technique we introduced reduces the truncation bias, making our results applicable to all sampling schemes. For a fixed $j$, Corollary \ref{cor-eigunif-fix} below discusses the uniform convergence rate under different ranges of $\bar{N}_2$. When $\alpha>5/2$, the truncation bias in equation \eqref{eq:eigunif} is dominated by the first two terms for all scenarios of $\bar{N}_2$. If $\bar{N}_2/(n/\log n)^{1/4}\rightarrow0$, $\|\hat\phi_{j}-\phi_{j}\|_{\infty}$ admits a typical one-dimensional kernel smoothing rate that  differs only by a $\log n$ factor, consistent with the result in \cite{hall2006propertieshans} and \cite{li2010uniform}. {In summary, the introduction of the double truncation technique facilitates the attainment of uniform convergence across all sampling schemes within functional data analysis. This advancement prompts more in-depth research, particularly in the analysis of non-compact data exhibiting in-curve dependence.}

\begin{corollary}\label{cor-eigunif-fix}
	Under the assumptions of Theorem \ref{thm-eigunif}. If $j$ is fixed and $\alpha>5/2$,
	\begin{itemize}
		\item[(1).] When $\bar{N}_2/(n/\log n)^{1/4}\rightarrow0 $ and  $h\asymp(n\bar{N}_2/\log n)^{-1/5} $,
		$$\sup_{t\in[0,1]}|\hat{\phi}_{j}(t)- {\phi}_{j}(t)| =O\left(  \left(\frac{\log n}{n\bar{N}_2} \right)^{\frac{2}{5}}\right) \quad \text{ a.s. }$$ 
		\item[(2).] When $\bar{N}_2/(n/\log n)^{1/4}\rightarrow C_2 >0 $ and $h\asymp(n/\log n)^{1/4} $
		$$ \sup_{t\in[0,1]}|\hat{\phi}_{j}(t)- {\phi}_{j}(t)| =O\left(  \sqrt{\frac{\log n}{n}}\right) \quad \text{ a.s. } $$
		\item[(3).] When $\bar{N}_2/(n/\log n)^{1/4}\rightarrow\infty$,  $h=o(n/\log n)^{1/4} $ and $\bar{N}_2h\rightarrow\infty$,
		$$ \sup_{t\in[0,1]}|\hat{\phi}_{j}(t)- {\phi}_{j}(t)| =O\left(  \sqrt{\frac{\log n}{n}}\right) \quad \text{ a.s. } $$
	\end{itemize}
\end{corollary}

The following corollary established the optimal uniform convergence rate for eigenfunctions with diverging indices under mild assumptions, {a finding that is new even in the fully observed scenario.} In contrast to the $\mathcal{L}^2$ convergence, the maximum number of eigenfunctions that can be well-estimated under the $|\cdot|_\infty$ norm is slightly smaller and depends on the moment {characterized by} assumption $\alpha$. {If $\alpha > {5}/{2}$, then $m_{\max}$ can increase  as the sample size $n$ goes to infinity, and the phase transition point $\bar{N}_2 \geq m_{\max}^a$ aligns with the $\mathcal{L}^2$ case.}

\begin{corollary}\label{cor-eigunif-diverg}
	Under the assumptions of Theorem \ref{thm-eigunif} and  \ref{asm:u1}. Given $\alpha>5/2$, denote $m_{\mathrm{max}}= \min \{n^{\frac{1}{2a}},n^{\frac{\alpha-5/2}{2\alpha a-\alpha-1}} \}$. If $\bar{N}_2\geq m_{\max}^a$, $h^{4}m_{\max}^{2c}\leq n^{-1}$ and  $\bar{N}_2h\geq m_{\max}$, for all $j\leq m_{\mathrm{max}}$, 
	$$ \sup_{t\in[0,1]}|\hat{\phi}_{j}(t)- {\phi}_{j}(t)|=O\left( \frac{j}{\sqrt{n}}(\sqrt{\ln n}+\ln j)\right)\quad a.s. $$
\end{corollary}

\section{Numerical experiment}\label{sec:simu}
In this section, we carry out a numerical evaluation of the convergence rates of eigenfunctions. The underlying trajectories are generated as $X_{i}(t) = \sum_{j=1}^{50} \xi_{ij} \phi_{j}(t), \text{ for } i = 1, \ldots, n$. The principal component scores  are independently generated following the distribution $\xi_{ij} \sim N(0, j^{-2})$ for all  $i$ and $j$. We define the eigenfunctions as $\phi_{1}(t) \equiv 1$, and for $j \geq 1$, as $\phi_{j}(t) = \sqrt{2}\cos(j\pi t)$. The actual observations are $X_{ij} = X_{i}(t_{ij}) + \epsilon_{ij}$, with noise $\epsilon_{ij}$ following $\mathcal{N}(0, 0.1^2)$, and the time points $t_{ij}$ sampled from a uniform distribution $\mathrm{Unif}[0, 1]$, for $j = 1, \ldots, N$. Each setting is repeated for 200 Mento-Carlo runs to mitigate the randomness that may occur in a single simulation.

\begin{figure}[h]
  \centering
  \includegraphics[width=0.55\linewidth]{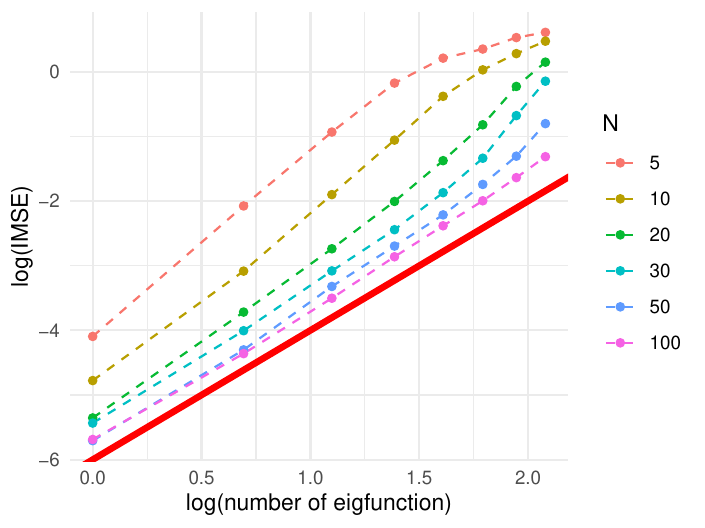}
  \caption{Plot of $\log(\|\hat\phi_{j}-\phi_j\|^2)$  over the logarithm of number of eigenfunctions. The sample size is $n=240$. The colored dashed lines correspond to different value of $N$. The solid red line represents the theoretical optimal value. }
  \label{fig:num}
\end{figure}

\begin{figure}[h]
  \centering
  \includegraphics[width=0.95\linewidth]{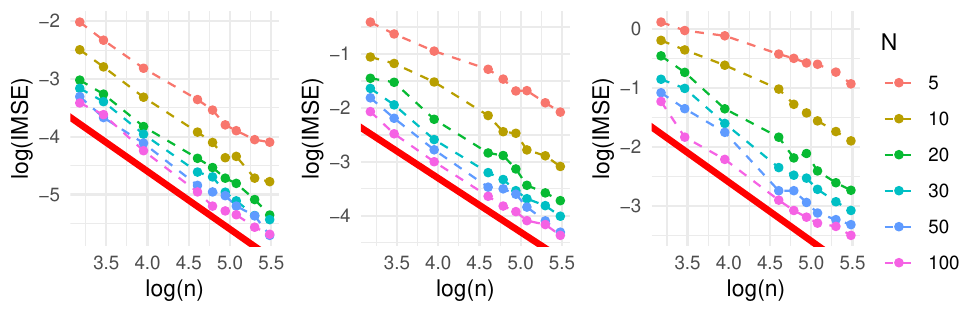}
  \caption{Plot of $\log(\|\hat\phi_{j}-\phi_j\|^2)$ for the first (left), second (middle) and third (right) eigenfunctions  over $log(n)$.  The colored dashed lines correspond to different value of $N$. The solid red line represents the theoretical optimal value. }
  \label{fig:eig1-3}
\end{figure}

When the phase transition occurs, our theory indicates a proportional relationship such that $\log(\|\hat\phi_{j}-\phi_j\|^2) \propto 2\log(j)$ for each fixed $n$ and $\log(\|\hat\phi_{j}-\phi_j\|^2) \propto -\log(n)$ for each fixed $j$. Figure \ref{fig:num} illustrates this phenomenon  by showing  that as $N$ increases, the relationship between $\log(\|\hat\phi_{j}-\phi_j\|^2)$ and $j$ tends to be linear with a slope of $2$, indicating that the phase transition might occur around $N=50$. Figure \ref{fig:num} additionally offers a practical way for identifying the phase transition point of \( N \). This is valuable in guiding both data collection and experimental design, contributing to more cost-effective data collecting  strategies.     Similarly, Figure \ref{fig:eig1-3}  shows that as $N$ increases, the relationship between $\log(\|\hat\phi_{j}-\phi_j\|^2)$ and $\log(n)$ also tends to follow a linear trend, but with a slope of $-1$.

\section{Conclusion and discussion}\label{sec:conc}
In this paper, we focus on the convergence rate of eigenfunctions with diverging indices for discretely observed functional data. We propose new techniques to handle the perturbation series and establish sharp bounds for eigenvalues and eigenfunctions across different convergence types. Additionally, we extend the partition ``dense'' and ``sparse'' defined for mean and covariance functions to principal components analysis. Another notable contribution of this paper is the double truncation technique for handling uniform convergence. 
%Due to the non-compactness of functional data, truncation is always necessary to apply the Bernstein inequality to achieve uniform convergence. 
Existing results on uniform convergence for covariance estimation require a strong moment condition on $X(t)$ and are only applicable to sparse functional data where $N_i=o_{p}(n^{1/4})$. By employing the double truncation technique proposed in this paper, we establish an improved bound for the truncated bias, which ensures the uniform convergence of the covariance and eigenfunctions across all sampling schemes under mild conditions. These asymptotic properties play a direct role in various types of statistical inference involving functional data \citep{yao2005jasa, li2010uniform}.

Furthermore, the optimal rate achieved in this paper holds significant implications for downstream analysis. Since most functional regression models encounter inverse issues due to the infinite dimensionality of functional covariates, the convergence rates in this paper would help improve existing theoretical findings in downstream analyses from fully observed functional data to various sampling designs.  Consider the functional linear model $\E[Y_{i}|X_{i}]=\int \beta X_{i}$ with $\beta=\sum_{k=1}^{\infty}k^{-b}\phi_{k}$. Without a sharp bound for eigenfunctions, achieving optimal convergence for standard and efficient plug-in estimators becomes challenging. Therefore, methods like approximated least squares and sample splitting, as discussed in \citep{zhou2023functional}, are necessary in the modeling phase. These complex methods require estimating principal components scores and do not efficiently utilize information gained from pooling. As a result, the phase transition for the functional linear model obtained by \cite{zhou2023functional} is ${\max\{n^{\frac{2a+2}{a+2b}} ,n^{\frac{a+b}{a+2b}} \}}$, which is significantly greater than $1/2$. In contrast, by using the new results in this paper, one can directly apply the plug-in estimator from \cite{hall2007methodology} and achieve a phase transition of $n^{\frac{1}{4}+\frac{3a-2b}{4(a+2b)}}$ for the functional linear model. Additionally, for complex regression models like the functional generalized linear model or functional Cox model, the methods developed in this paper could serve as a cornerstone for further exploration.

\begin{appendix}

	\section*{Proof of Theorem \ref{thm-eig}}
Recall $\hat{C}(s,t)=\hat{\beta}_{0}$ in the optimization problem \eqref{eq:ll} and $\hat{\beta}_0$ has the analytical from
	$$\hat{\beta}_{0}(s,t)=\frac{R_{00}(s,t)I_{1}(s,t)+R_{10}(s,t)I_2(s,t)+R_{01}(s,t)I_3(s,t)}{S_{00}(s,t)I_{1}(s,t)+S_{10}(s,t)I_2(s,t)+S_{01}(s,t)I_3(s,t)},$$
where
\begin{align*}
	S_{pq}(s,t)=&\sum_{i=1}^{n}v_{i}\sum_{l_{1}\neq l_{2}}^{N_{i}}\frac{1}{h^2}\K\left(\frac{t_{il_1}-s}{h} \right) \K\left(\frac{t_{il_2}-t}{h} \right) \left(\frac{t_{il_1}-s}{h} \right)^p \left(\frac{t_{il_2}-t}{h} \right )^{q}\\
	R_{pq}(s,t)=&\sum_{i=1}^{n}v_{i}\sum_{l_{1}\neq l_{2}}^{N_{i}}\frac{1}{h^2}\K\left(\frac{t_{il_1}-s}{h} \right) \K\left(\frac{t_{il_2}-t}{h} \right)\left(\frac{t_{il_1}-s}{h} \right)^p \left(\frac{t_{il_2}-t}{h} \right )^{q}\delta_{il_1l_2}
\end{align*}
for $p,q=0,1,2$ and
\begin{align*}
	&I_{1}(s,t)=S_{02}(s,t)S_{20}(s,t)-S_{11}^2(s,t),\ I_{2}(s,t)=S_{01}(s,t)S_{11}(s,t)-S_{02}(s,t)S_{10}(s,t)\\
	&I_{3}(s,t)=S_{10}(s,t)S_{11}(s,t)-S_{20}(s,t)S_{01}(s,t).
\end{align*} 
Some calculations show that
$$\hat{\beta}_{1}(s,t)=\frac{1}{h}\frac{R_{00}(s,t)I_{2}(s,t)+R_{10}(s,t)J_1(s,t)+R_{01}(s,t)J_3(s,t)}{S_{00}(s,t)I_{1}(s,t)+S_{10}(s,t)I_2(s,t)+S_{01}(s,t)I_3(s,t)}$$
and 
$$\hat{\beta}_{2}(s,t)=\frac{1}{h}\frac{R_{00}(s,t)I_{3}(s,t)+R_{10}(s,t)J_3(s,t)+R_{01}(s,t)J_2(s,t)}{S_{00}(s,t)I_{1}(s,t)+S_{10}(s,t)I_2(s,t)+S_{01}(s,t)I_3(s,t)}$$
where 
\begin{align*}
	&J_{1}(s,t)=S_{00}(s,t)S_{02}(s,t)-S_{01}^2(s,t),\ J_{2}(s,t)=S_{00}(s,t)S_{20}(s,t)-S_{10}^2(s,t)\\
	&J_{3}(s,t)=S_{10}(s,t)S_{01}(s,t)-S_{00}(s,t)S_{11}(s,t).
\end{align*}	
In the following, we omit the arguments $(s,t)$ in the functions  $\hat{\beta}_2$, $\hat{\beta}_3$, $S_{pq},\ R_{pq}$, $I_r$ and $J_r$ for $p,q=0,1,2$ and $r=1,2,3$ when there is no ambiguity. Simple calculations show that $\hat{C}(s,t)=(R_{00}-h\hat{\beta}_1S_{10}-h\hat{\beta}_2S_{01})/S_{00}$. We further denote 
$$\tilde{C}_{0}(s,t)=\left\{R_{00}-C(s,t)S_{00}-h\frac{\partial C(s,t)}{\partial s }(s,t)S_{10}- h\frac{\partial C(s,t)}{\partial t }(s,t)S_{01}\right\}/{f(s)f(t)},  $$
where $f(s)$ is the density function of $t_{ij}$ in the case of the random design, and $f(s)=1_{\{0\leq s\leq 1\}}$ in the case of the fixed regular design. The subsequent proof is structured into three steps to ensure clearer understanding.

\textbf{Step 1: Discrepancy between $\hat{C}(s,t)$ and $\tilde{C}_{0}$.}
The following lemma bound the discrepancy between $\hat{C}(s,t)-C(s,t)$ and $\tilde{C}_{0}(s,t)$, and its proof can be found in Section S3 of the supplement. 
\begin{lemma}\label{lem:dis-C}
	Under Assumption \ref{asm:a1}, \ref{asm:cov} and \ref{asm:u1}.
	\begin{itemize}
		\item[(a)] For the random design, if Assumption \ref{asm:r} holds,
		\begin{equation}\label{eq:ll-r1}
	\left \|\hat{C}(s,t)-C(s,t) -\tilde{C}_{0}(s,t) \right\rhsnorm ^2=O_{P}\left( h^6+{\frac{h^2}{n}\left\{1+\frac{1}{\bar{N}^2_{2}h^2} \right\} } \right).
\end{equation}
		\item[(b)] For the fix design, if Assumption \ref{asm:e} and \ref{asm:e-c} hold,
		\begin{equation}\label{eq:ll-e1}
	\left \|\hat{C}(s,t)-C(s,t) -\tilde{C}_{0}(s,t) \right\rhsnorm ^2=O_{P}\left( h^6+{\frac{h^2}{n}\left\{1+\frac{1}{N^2h^2} \right\} } \right).
\end{equation}
	\end{itemize}
\end{lemma}

\textbf{Step 2: bound the projection $\E( [\iint\tilde{C}_{0}(s,t)\phi_{j}(s)\phi_{k}(t) \diff s\diff t   ]^2  )$. }
In this part, we will prove the following lemma, which  provides expectation bounds for the projections of $\tilde{C}_{0}(s,t)$ with respect to $\phi_j$ and $\phi_k$.
\begin{lemma}\label{lem:Cjk}
	Under assumptions \ref{asm:a1} to \ref{asm:a3}, $h^{4}j^{2a+2c}\lesssim1$ and $hj^{a}\log n\lesssim1 $,  
	\begin{itemize}
		\item[(a)] If Assumption \ref{asm:r} holds, for all $1\leq k\leq 2j$,
			\begin{align*}
				&\E\left[ \left\{\iint\tilde{C}_{0}(s,t)\phi_{j}(s)\phi_{k}(t) \diff s\diff t   \right \}^2  \right ]\lesssim  \frac{1}{n}\left(j^{-a}k^{-a}+\frac{j^{-a}+k^{-a}}{\bar{N}_2}+\frac{1}{\bar{N}^2_2} \right)+h^{4}k^{2c-2a}
			\end{align*}
	and
	\begin{align*}
		&\sum_{k=j+1}^{\infty}\E\left[ \left\{\iint\tilde{C}_{0}(s,t)\phi_{j}(s)\phi_{k}(t) \diff s\diff t   \right \}^2  \right ]\\
		\lesssim& \frac{1}{n}\left(j^{1-2a}+\frac{h^{-1}j^{-a}+j^{1-a}}{\bar{N}_2^2}+\frac{1}{h\bar{N}^2_2} \right)+h^{4}j^{1+2c-2a}.
	\end{align*}
	\item[(b)] If Assumptions \ref{asm:e} and \ref{asm:e-c} hold, for all $1\leq k\leq 2j$,
	\begin{align*}
				&\E\left[ \left\{\iint\tilde{C}_{0}(s,t)\phi_{j}(s)\phi_{k}(t) \diff s\diff t   \right \}^2  \right ]\lesssim  \frac{1}{n}\left(j^{-a}k^{-a}+\frac{j^{-a}+k^{-a}}{N}+\frac{1}{N^{2}} \right)+h^{4}k^{2c-2a}
			\end{align*}
	and
	\begin{align*}
		&\sum_{k=j+1}^{\infty}\E\left[ \left\{\iint\tilde{C}_{0}(s,t)\phi_{j}(s)\phi_{k}(t) \diff s\diff t   \right \}^2  \right ]\\
		\lesssim& \frac{1}{n}\left(j^{1-2a}+\frac{h^{-1}j^{-a}+j^{1-a}}{N}+\frac{1}{hN^{2}} \right)+h^{4}j^{1+2c-2a}.
	\end{align*}
	\end{itemize}
\end{lemma}

\begin{proof}[Proof of Lemma \ref{lem:Cjk}]
We focus on the first statement of Lemma \ref{lem:Cjk}. The proof of the fix regular design is similar and we put it in supplement. By definition of $\tilde{C}_{0}(s,t)$, we need to bound the bias and variance terms of 
\begin{equation}\label{eq:cjk2}
	\iint\left\{R_{00}-C(s,t)S_{00}-h\frac{\partial C}{\partial s}(s,t)S_{10}-h\frac{\partial C}{\partial t}(s,t) S_{01} \right\}\frac{\phi_{j}(s)\phi_{k}(t)}{f(s)f(t)} \diff s\diff t.
\end{equation}
For the random design case, by analogous calculation as proof of Theorem 3.2 in \cite{zhang2016sparse}, one has
\begin{align*}
	&\E\left[R_{00}-C(s,t)S_{00}-h\frac{\partial C}{\partial s}(s,t)S_{10}-h\frac{\partial C}{\partial t}(s,t) S_{01}\right]\\
	=&\frac{h^2}{2}\mathrm{K}_{2} \frac{\partial C(s,t)}{\partial s^2}(s,t)f(s)f(t)+\frac{h^2}{2}\mathrm{K}_{2} \frac{\partial C(s,t)}{\partial t^2}(s,t)f(s)f(t)+O(h^3),
\end{align*} 
where $\K_2=\int u^2 \K(u)\diff u$. Thus, for the bias part of equation \eqref{eq:cjk2}, for all $k\leq2j$,
\begin{equation}\label{eq:cjk3}
\begin{aligned}
	&\left(\E \iint\left\{R_{00}-C(s,t)S_{00}-h\frac{\partial C}{\partial s}(s,t)S_{10}-h\frac{\partial C}{\partial t}(s,t) S_{01} \right\}\frac{\phi_{j}(s)\phi_{k}(t)}{f(s)f(t)} \diff s\diff t\right)^2\\
%	=& \frac{h^{4}}{4}\K_{2}^2\left[ \iint\left\{ \frac{\partial^2 C}{\partial t^2}(t,s)+\frac{\partial^2 C}{\partial t^2}(s,t)\right\}\phi_{j}(s)\phi_{k}(t)\diff s\diff t\right]^2+o(j^{-2a}h^4) \\
	=& \frac{h^{4}}{2}\K_{2}^2 \left[\iint \sum_{r=1}^{\infty}\lambda_{r}\phi_{r}(s)\phi_{r}^{(2)}(t)\phi_{j}(s)\phi_{k}(t)\diff s\diff t\right]^2+o(j^{-2a}h^4) \\
	\leq&\frac{h^{4}}{2}\K_{2}^2 \lambda_{j}^2\|\phi_{j}\|^2_{\infty}+o(h^4) 
	=O( h^{4}k^{-2a}j^{2c}).
\end{aligned}
\end{equation}

For the tail summation, similarly 
\begin{equation}\label{eq:cjk4}
	\begin{aligned}
		&\sum_{k\geq j}\left(\E \iint\left\{R_{00}-CS_{00}-h\frac{\partial C}{\partial s}(s,t)S_{10}-h\frac{\partial C}{\partial t}(s,t) S_{01} \right\}\frac{\phi_{j}(s)\phi_{k}(t)}{f(s)f(t)} \diff s\diff t\right)^2\\
%		=& \frac{h^{4}}{2}\K_{2}^2\sum_{k>j}^{\infty} \left[ \iint \frac{\partial^2 C}{\partial t^2}(t,s)\phi_{j}(s)\phi_{k}(t)\diff s\diff t\right]^2+o(j^{-2a}h^4) \\
		\leq& \frac{h^{4}}{2}\K_{2}^2\left\|\int \frac{\partial^2 C}{\partial t^2}(t,s)\phi_{j}(s)\diff s \right\rhsnorm+o(j^{-2a}h^4)\\
%		=& \frac{h^{4}}{2}\K_{2}^2\left\|\int \sum_{r=1}^{\infty}\lambda_{r}\phi_{r}(s)\phi_{r}^{(2)}(t) \phi_{j}(s)\diff s \right\rhsnorm+o(j^{-2a}h^4)\\
		=&O( h^{4}j^{1+2c-2a}).
	\end{aligned}
\end{equation}

For the variance of equation \eqref{eq:cjk2}, note that
\begin{equation}\label{eq:cjk5}
	\begin{aligned}
		&\var \left(  \iint\left\{R_{00}-CS_{00}-h\frac{\partial C}{\partial s}(s,t)S_{10}-h\frac{\partial C}{\partial t}(s,t) S_{01} \right\}\frac{\phi_{j}(s)\phi_{k}(t)}{f(s)f(t)} \diff s\diff t \right)\\
		\leq&\E \left[ \left\{\iint R_{00}(s,t)\frac{\phi_{j}(s)\phi_{k}(t)}{f(s)f(t)} \diff s\diff t \right\}^2 \right]+\E \left[ \left\{\iint C(s,t)S_{00}(s,t)\frac{\phi_{j}(s)\phi_{k}(t)}{f(s)f(t)} \diff s\diff t \right\}^2 \right]\\
		&+h^{2}\E \left[ \left\{\iint \frac{\partial C}{\partial s}(s,t)S_{10}(s,t)\frac{\phi_{j}(s)\phi_{k}(t)}{f(s)f(t)} \diff s\diff t \right\}^2 \right]\\
		&+h^{2}\E \left[ \left\{\iint \frac{\partial C}{\partial t }(s,t)S_{01}(s,t)\frac{\phi_{j}(s)\phi_{k}(t)}{f(s)f(t)} \diff s\diff t \right\}^2 \right].
	\end{aligned}
\end{equation}
We start with the first term in the right hand side of equation \eqref{eq:cjk5}. To simplify the notation, we shall introduce the following notation:
$$\mathcal{T}_{h}f(x)=\frac{1}{h}\int \K\left(\frac{x-y}{h} \right)f(y)\diff y. $$
Then
\begin{align*}
	&\E \left[ \left\{\iint R_{00}(s,t)\frac{\phi_{j}(s)\phi_{k}(t)}{f(s)f(t)} \diff s\diff t \right\}^2 \right]
	=\E\left[\left\{  \sum_{i=1}^{n}v_{i}\sum_{l_1\neq l_2}^{N_i}\delta_{il_1l_2}\mathcal{T}_{h}\frac{\phi_{j}}{f}(t_{il_1})\mathcal{T}_{h}\frac{\phi_{k}}{f}(t_{il_2}) \right\}^2\right]\\
	=&\sum_{i=1}^{n}v_{i}^2\left\{ 4!\binom{N_{i}}{4}A_{i1}+3!\binom{N_{i}}{3}A_{i2}+2!\binom{N_{i}}{2}A_{i3} \right\}
\end{align*}
with
\begin{align*}
	A_{i1}=&\mathbb{E}\left(\left\langle X_{i}f,\mathcal{T}_{h}\frac{\phi_{j}}{f}\right \rangle^2\left\langle X_{i}f,\mathcal{T}_{h}\frac{\phi_{k}}{f}\right\rangle^2 \right ),\\
	A_{i2}=&2\mathbb{E}\left\{\left(\left\langle X_{i}f\mathcal{T}_{h}\frac{\phi_{j}}{f},X_{i}\mathcal{T}_{h}\frac{\phi_{k}}{f}\right \rangle +\sigma_{X}^2\left\langle  \mathcal{T}_{h}\frac{\phi_{j}}{f},f\mathcal{T}_{h}\frac{\phi_{k}}{f}\right\rangle \right)\right.\left.\left\langle X_{i}f,\mathcal{T}_{h}\frac{\phi_{j}}{f}\right\rangle\left\langle X_{i}f,\mathcal{T}_{h}\frac{\phi_{k}}{f}\right\rangle\right\}\\
		&+\mathbb{E}\left\{\left\langle X_{i}f, \mathcal{T}_{h}\frac{\phi_{j}}{f}\right\rangle^2\left\langle (X_{i}^{2}+\sigma_{X}^2)f, \left(\mathcal{T}_{h}\frac{\phi_{k}}{f} \right)^2\right\rangle\right\}\\&+\mathbb{E}\left\{\left \langle X_{i}f, \mathcal{T}_{h}\frac{\phi_{k}}{f}\right\rangle^2\left\langle (X_{i}^{2}+\sigma_{X}^2)f, \left(\mathcal{T}_{h}\frac{\phi_{j}}{f} \right)^2\right\rangle\right\}\\
		:=&A_{i21}+A_{i22}+A_{i23}\\
	A_{i3}=&\mathbb{E}\left\{ \left\langle (X_i^2+\sigma_{X}^2)f,   \left(\mathcal{T}_{h}\frac{\phi_{j}}{f} \right)^2\right\rangle \left\langle (X_i^2+\sigma_{X}^2)f,   \left(\mathcal{T}_{h}\frac{\phi_{k}}{f} \right)^2\right\rangle \right\} \\
			&+\mathbb{E} \left (\left\langle X_{i}\mathcal{T}_{h}\frac{\phi_{j}}{f},X_{i}f\mathcal{T}_{h}\frac{\phi_{k}}{f}\right\rangle+\sigma_{X}^{2}\left \langle  \mathcal{T}_{h}\frac{\phi_{j}}{f}, f\mathcal{T}_{h}\frac{\phi_{k}}{f}\right \rangle  \right)^2 \\
			:=&A_{i31}+A_{i32}.
\end{align*}
By Cauchy--Schwarz and AM--GM inequality, 
\begin{align*}
		A_{i21}
		\leq&\mathbb{E}\left\{\left \langle X_if, \mathcal{T}_{h}\frac{\phi_{j}}{f}\right \rangle^2\left\langle (X_{i}^{2}+\sigma_{X}^2)f, \left(\mathcal{T}_{h}\frac{\phi_{k}}{f} \right)^2\right\rangle\right\}\\&+\E\left\{\left\langle X_if, \mathcal{T}_{h}\frac{\phi_{k}}{f}\right\rangle^2\left\langle (X_{i}^{2}+\sigma_{X}^2)f, \left(\mathcal{T}_{h}\frac{\phi_{j}}{f} \right)^2\right\rangle\right\}\\
		=&A_{i22}+A_{i23}.	
\end{align*}
Similarly $A_{i32}\leq A_{i31}$, thus, $A_{i2}\leq 2\{A_{i22}+A_{i23} \},$ $ A_{i3}\leq 2A_{i31}$. Combine all above,
\begin{equation}\label{eq:cjk6}
	\begin{aligned}
		&\E \left[ \left\{\iint R_{00}(s,t)\frac{\phi_{j}(s)\phi_{k}(t)}{f(s)f(t)} \diff s\diff t \right\}^2 \right]\\
		\lesssim&\sum_{i=1}^{n}v_{i}^2\left\{4! \binom{N_{i}}{4}A_{i1}+3!\binom{N_{i}}{3}(A_{i22}+A_{i23})+2!\binom{N_{i}}{2}A_{i31}\right\}.
	\end{aligned}
\end{equation}
The following lemma in bounding  $A_{i1},A_{i22},A_{i23}$ and $A_{i31}$, its proof can be found in the supplement. 
\begin{lemma}\label{lem-4thmoment}
	Under assumptions \ref{asm:a1} to \ref{asm:a3} and {$h^{4}j^{2c+2a}\lesssim1$}, $hj^{a}\log n\lesssim1 $, there is   $\mathbb{E}(\langle X_if,\mathcal{T}_h\frac{\phi_{k}}{f} \rangle^4)\lesssim k^{-2a}\ \text{for}\ 1\leq k\leq 2j$ and
	$$
\mathbb{E}\left(\sum_{k>j}\left\langle X_if,\mathcal{T}_h\frac{\phi_k}{f} \right \rangle^2\right)^2 \lesssim j^{2-2a}.$$
\end{lemma}

By Lemma \ref{lem-4thmoment} and Cauchy-Schwarz inequality
\begin{equation}\label{eq:varA1}
\begin{aligned}
	&A_{i1} \leq \left(\mathbb{E}\left\langle X_if,\mathcal{T}_{h}\frac{\phi_j}{f}\right\rangle^4\mathbb{E}\left\langle X_if,\mathcal{T}_{h}\frac{\phi_k}{f}\right\rangle^4\right)^{{1}/{2}}\lesssim j^{-a}k^{-a}\text{ and }\\
	&\sum_{k>j} A_{i1}\leq \left(\mathbb{E}\left\langle X_if,\mathcal{T}_{h}\frac{\phi_j}{f}\right \rangle^4\mathbb{E}\left(\sum_{k>j}\left\langle X_if,\mathcal{T}_{h}\frac{\phi_k}{f}\right \rangle^2\right)^2\right)^{{1}/{2}}\lesssim  j^{1-2a}.
\end{aligned}
\end{equation}

For $A_{i22} $, by Lemma \ref{lem-4thmoment} and Cauchy-Schwarz inequality,
\begin{equation}\label{eq:varA22-1}
	\begin{aligned}
A_{i22}
\leq& \left\|\frac{\phi_{k}}{f}\right\|_{{\infty}}^2 \|f\|_{\infty}\left\{\mathbb{E}\left\langle X_if,\mathcal{T}_{h}\frac{\phi_j}{f}\right\rangle^4\mathbb{E}(\|X_i\|^2+\sigma_{X}^2)^2\right\}^{1/2}\lesssim  j^{-a},\\
\end{aligned}
\end{equation}
For the summation $\sum_{k>j}A_{i22}$, by Cauchy-Schwarz inequality
\begin{align*}
		&\sum_{k=1}^{\infty}\left|\mathcal{T}_{h}\frac{\phi_k}{f}\right|^2\leq\frac{1}{h^2}\int \left|\mathrm{K}\left(\frac{x-y}{h}\right)\right|^2\frac{1}{f^2(y)} \diff y\lesssim {h}^{-1} \text{ for all }x\in [0,1] .
	\end{align*}
Thus,
\begin{equation}\label{eq:varA22-2}
	\begin{aligned}
		&\sum_{k=j}^{\infty}A_{i22}=\mathbb{E}\left\{\left \langle X_if,\mathcal{T}_{h}\frac{\phi_j}{f}\right \rangle^2\sum_{k=j}^{\infty}\left\langle (X_{i}^{2}+\sigma_{X}^2)f, \left(\mathcal{T}_{h}\frac{\phi_{k}}{f} \right)^2\right\rangle\right\}\\
		\leq&\|f\|_{\infty}\sum_{k=1}^{\infty}\left\|\mathcal{T}_{h}\frac{\phi_k}{f}\right \|^2 \mathbb{E}\left\{\left\langle X_if,\mathcal{T}_{h}\frac{\phi_j}{f}\right \rangle^2(\|X_i\|^2+\sigma_{X}^2)\right\}\\
\lesssim&  h^{-1}\left(\mathbb{E}\left\langle X_if,\mathcal{T}_{h}\frac{\phi_j}{f}\right \rangle^4\right)^{1/2}\lesssim h^{-1}j^{-a}.
	\end{aligned}
\end{equation}

For $A_{i23} $, by Lemma \ref{lem-4thmoment}, 
\begin{equation}\label{eq:varA23-1}
	\begin{aligned} &A_{i23}=A_{i22}\leq \left\|\mathcal{T}_{h}\frac{\phi_j}{f}\right\|_{{\infty}}^2\|f\|_{\infty} \mathbb{E}\left\{\left\langle X_if,\mathcal{T}_{h}\frac{\phi_k}{f}\right \rangle^2(\|X_i\|^2+\sigma_{X}^2)\right\}\\\lesssim&\left(\mathbb{E}\left\langle X_if,\mathcal{T}_{h}\frac{\phi_k}{f}\right \rangle^4\right)^{1/2}\lesssim k^{-a},\quad \forall \ 1\leq k\leq 2j,
\end{aligned}
\end{equation}
and
\begin{equation}\label{eq:varA23-2}
	\begin{aligned} &\sum_{k>j}A_{i23}\leq\left\|\mathcal{T}_{h}\frac{\phi_j}{f}\right \|_{{\infty}}^2\|f\|_{\infty} \mathbb{E}\left[\sum_{k>j}\!\left\langle X_if,\mathcal{T}_{h}\frac{\phi_k}{f}\right\rangle^2(\|X_i\|^2+\sigma_{X}^2)\right]\\\lesssim& \left\{\mathbb{E}\left(\sum_{k>j}\!\left\langle X_if,\mathcal{T}_{h}\frac{\phi_k}{f}\right \rangle^2\right)^2\E(\|X\|^{2}+\sigma_{X}^2)^{2} \right\}^{\frac{1}{2}}\lesssim  j^{1-a}.
\end{aligned}
\end{equation}

For the last term $A_{i31} $, note that
 \begin{align*} A_{i31}=&\mathbb{E}\left\{ \left\langle (X_i^2+\sigma_{X}^2)f,   \left(\mathcal{T}_{h}\frac{\phi_{j}}{f} \right)^2\right\rangle \left\langle (X_i^2+\sigma_{X}^2)f,   \left(\mathcal{T}_{h}\frac{\phi_{k}}{f} \right)^2\right\rangle \right\}=O(1),
\end{align*}
and
\begin{equation}\label{eq:varA31-2}
	\begin{aligned}
	&\sum_{k=j}^{\infty}A_{i31}
\leq \mathbb{E} \left\{\left\langle (X_i^2+\sigma_{X}^2)f,   \left(\mathcal{T}_{h}\frac{\phi_{j}}{f} \right)^2\right\rangle (\left\|X_i\right \|^2\|f\|_{\infty}+\sigma_{X}^2\|f\|^2)\frac{\|\mathrm{K}\|^2}{h}\right\}\\
\lesssim& h^{-1}\mathbb{E}(\|X_i\|^2\|f\|_{\infty}+\sigma_{X}^2\|f\|_{\infty})^2\lesssim  h^{-1}.
\end{aligned}
\end{equation}

Combine equation \eqref{eq:varA1} to \eqref{eq:varA31-2}, for all $k\leq 2j$
\begin{align*}
	&\E \left[ \left\{\iint R_{00}(s,t)\frac{\phi_{j}(s)\phi_{k}(t)}{f(s)f(t)} \diff s\diff t \right\}^2 \right]
	\lesssim\frac{1}{n}\left( {j^{-a}k^{-a}}+\frac{j^{-a}+k^{-a}}{\bar{N}_{2}}+\frac{1}{\bar{N}_2^2} \right)
\end{align*}
and
\begin{align*}
	&\sum_{k>j}\E \left[ \left\{\iint R_{00}(s,t)\frac{\phi_{j}(s)\phi_{k}(t)}{f(s)f(t)} \diff s\diff t \right\}^2 \right]
	\lesssim\frac{1}{n}\left( {j^{1-2a}}+\frac{j^{-a}h^{-1}+j^{1-a}}{\bar{N}_{2}}+\frac{h^{-1}}{\bar{N}^2_2} \right).
\end{align*}
By similar analysis, the second term in the right hand side of equation \eqref{eq:cjk5} has the same convergence rate as the first term. Under $h^{4}j^{2a+2c}\lesssim1$ and  $hj^{a}\log n\lesssim1$, the last two terms in the right hand side of equation \eqref{eq:cjk5} are dominated by the first two terms. Then the proof of the random design case is complete by
\begin{align*}
	&\E \left(\iint\left\{R_{00}-C(s,t)S_{00}-h\frac{\partial C}{\partial s}(s,t)S_{10}-h\frac{\partial C}{\partial t}(s,t) S_{01} \right\}\frac{\phi_{j}(s)\phi_{k}(t)}{f(s)f(t)} \diff s\diff t \right)^2\\
	\lesssim& \frac{1}{n}\left( {j^{-a}k^{-a}}+\frac{j^{-a}+k^{-a}}{\bar{N}_{2}}+\frac{1}{\bar{N}^2_2} \right)+h^{4}k^{-2a}j^{2c}
\end{align*}
for all $k\leq 2j$ and 
\begin{align*}
	&\sum_{k>j}\E \left(\iint\left\{R_{00}-C(s,t)S_{00}-h\frac{\partial C}{\partial s}(s,t)S_{10}-h\frac{\partial C}{\partial t}(s,t) S_{01} \right\}\frac{\phi_{j}(s)\phi_{k}(t)}{f(s)f(t)} \diff s\diff t \right)^2\\
	\lesssim& \frac{1}{n}\left( {j^{1-2a}}+\frac{j^{-a}h^{-1}+j^{1-a}}{\bar{N}_{2}}+\frac{h^{-1}}{\bar{N}^2_2} \right)+h^{4}j^{1-2a+2c}.
\end{align*}
\end{proof}

\textbf{Step 3: perturbation series.}
	By the proof of Theorem 5.1.8 in \cite{hsing2015theoretical}, for $j\in m $ and on the set $\Omega_{m}(n,N,h)=\{\|\hat{C}-C \rhsnorm\leq \eta_{m}/2 \}$, we have the following expansion,
\begin{equation}\label{eq:eig-1}
	\begin{aligned}
		\hat{\phi}_{ j}-\phi_{j}
	= & \sum_{k\neq j}\frac{ \int(\hat{C}  -C) \phi_{j} \phi_{k}}{\left(\lambda_{j}-\lambda_{k}\right)}\phi_{k}+\sum_{k\neq j}\frac{\int(\hat{C}  -C)(\hat\phi_{ j}- \phi_{j}) \phi_{k}}{\left(\lambda_{j}-\lambda_{k}\right)}\phi_{k} \\
	&+\sum_{k\neq j}\sum_{s=1}^{\infty}\frac{(\lambda_{j}-\hat\lambda_{ j})^s}{(\lambda_{j}-\lambda_{k} )^{s+1} }\left\{\int(\hat{C}  -C) \hat\phi_{ j} \phi_{k}\right\} \phi_{k}+\left\{\int (\hat\phi_{ j}- \phi_{j}) \phi_{j}\right\}\phi_{j}.
	\end{aligned}
\end{equation}
Such kind of expansion can also be found in \cite{hall2006properties} and \cite{li2010uniform}. Below, when we say that a bound is valid when $\Omega_{m}(n,N,h)$ holds, this should be interpreted as stating that the bound is valid for all realizations for which $\|\hat{C}-C \rhsnorm\leq \eta_{m}/2 $ \citep{hall2007methodology}. Under assumptions $m^{2a+2}/n\rightarrow0$, $m^{2a+2}/(n\bar{N}_2^2h^2)\rightarrow0 $ and $h^{4}\max\{ m^{2a+2c},m^{4a}\log n \}\lesssim1 $, we have $\mathbb{P}(\Omega_{m}(n,N,h))\rightarrow1$. Since $\lim_{D\rightarrow\infty}\limsup_{n\rightarrow\infty} P(\|\hat{\phi}_j-\phi_j\|>D\tau_{n} )=0 $ implies $\|\hat{\phi}_j-\phi_j\|=O_{P}(\tau_{n}) $, thus the results in $O_{P}$ form we want to prove relate only to probabilities of differences. It suffices to work with bounds that are established under the assumption such that  $\mathbb{P}(\Omega_{m})\rightarrow1$ holds. \citep[Section 5.1]{hall2007methodology}. 

 We first show that $\E(\|\hat\phi_{j}-\phi_{j}\|^{2})$ is dominated by the $\mathcal{L}^{2}$ norm of the first term in the right hand side of equation \eqref{eq:eig-1}. By Bessel's inequality, we see that 
\begin{equation}\label{eq:eig-2}
	\mathbb{E}\left\| \sum_{k\neq j}\frac{\int(\hat{C}  -C)(\hat\phi_{ j}- \phi_{j}) \phi_{k}}{\left(\lambda_{j}-\lambda_{k}\right)}\phi_{k} \right\|^2\leq \mathbb{E}\frac{\|\hat{C}  -C\rhsnorm^2\|\hat\phi_{ j}-\phi_{j} \|^2}{(2\eta_{j})^2}<\frac{1}{16}\mathbb{E}\|\hat\phi_{ j}-\phi_{j} \|^2,
\end{equation}
where the last equality comes from the fact $\eta_{j}^{-1}\|\hat{C}  -C\|<1/2$ on $\Omega_{m}(n,N,h)$. Similarly,
{\small{\begin{equation}\label{eq:eig-3}
	\begin{aligned}
		&\mathbb{E}\left\|\sum_{k\neq j}\sum_{s=1}^{\infty}\frac{(\lambda_{j}-\hat\lambda_{ j})^s}{(\lambda_{j}-\lambda_{k} )^{s+1} }\left\{\int(\hat{C}  -C) \hat\phi_{ j} \phi_{k}\right\} \phi_{k} \right\|^2\\
	 	=&\mathbb{E}\sum_{k\neq j}\frac{(\lambda_{j}-\hat\lambda_{ j} )^2}{(\lambda_{j}-\lambda_{k} )^2(\hat\lambda_{ j}-\lambda_{k} )^2}\left\{\int(\hat{C}  -C) \hat\phi_{ j} \phi_{k}\right\}^2\\
	 	\leq &2\mathbb{E}\frac{\|\hat{C}  -C\rhsnorm^2}{(2\eta_{j}-\|\hat{C}  -C\rhsnorm)^2}\left[\sum_{k\neq j}\frac{\left\{\int(\hat{C}  -C) \phi_{j} \phi_{k}\right\}^{2} }{\left(\lambda_{j}-\lambda_{k}\right)^{2}}\right.+\left.\sum_{k\neq j}\frac{\left\{\int(\hat{C}  -C)(\hat\phi_{ j} -\phi_{j}) \phi_{k}\right\}^{2} }{\left(\lambda_{j}-\lambda_{k}\right)^{2}} \right]\\
	 	\leq & \frac{8}{9}\mathbb{E}\left[ \frac{\|\hat{C}  -C\rhsnorm^2}{\eta_{j}^2} \sum_{k\neq j}\frac{\left\{\int(\hat{C}  -C) \phi_{j} \phi_{k}\right\}^{2} }{\left(\lambda_{j}-\lambda_{k}\right)^{2}}+ \frac{\|\hat{C}  -C\rhsnorm^4}{\eta_{j}^4}\|\hat\phi_{j}-\phi_{j} \|^2\right]\\
	 	\leq & \frac{2}{9}\mathbb{E}\sum_{k\neq j}\frac{\left\{\int(\hat{C}  -C) \phi_{j} \phi_{k}\right\}^{2} }{\left(\lambda_{j}-\lambda_{k}\right)^{2}}+\frac{1}{18}\mathbb{E}\|\hat\phi_{j}-\phi_{j} \|^2 .	 \end{aligned}
\end{equation}}}
Combing \eqref{eq:eig-1} to \eqref{eq:eig-3} and the fact $\|\{\int (\hat\phi_{j}- \phi_{j}) \phi_{j}\}\phi_{j}\|=1/2 \|\hat\phi_{j}-\phi_{j} \|^2 $, $\mathbb{E}(\|\hat\phi_{j}-\phi_{j}\|^2)$ is dominated by the first term in the right hand  side of equation \eqref{eq:eig-1}. Thus,
\begin{equation}\label{eq:eig-31}
\begin{aligned}
	&\mathbb{E}(\|\hat\phi_{j}-\phi_{j}\|^2)
	\lesssim& \sum_{k\neq j}\frac{\mathbb{E}\left[\left\{\int(\hat{C}  -C - \tilde{C}_{0}) \phi_{j} \phi_{k}\right\}^{2}\right] }{(\lambda_{j}-\lambda_{k})^2}+\sum_{k\neq j}\frac{\mathbb{E}\left\{\left(\int\tilde{C}_{0} \phi_{j} \phi_{k}\right)^{2}\right\} }{(\lambda_{j}-\lambda_{k})^2}.
\end{aligned}
\end{equation}

For the random design,  under Assumption \ref{asm:r} and {$ hj^{a}\log n\lesssim1$}, by Bessel equality and  \eqref{eq:ll-r1},  the first term in the right hand side of equation \eqref{eq:eig-31} is bounded by 
\begin{align*}
	&\sum_{k\neq j}\frac{\mathbb{E}\left\{\int(\hat{C}  -C - \tilde{C}_{0}) \phi_{j} \phi_{k}\right\}^{2} }{(\lambda_{j}-\lambda_{k})^2}\leq \frac{\left\| \hat{C}  -C - \tilde{C}_{0}\right\rhsnorm^2 }{\eta_{j}^2}\\
	=&O_{P}\left(\frac{j^{2}}{n}\left\{1+\frac{j^a}{\bar{N}_2}+\frac{j^{2a}}{\bar{N}^2_2}  \right\} \right)+o_{P}\left(h^{2}j^{2c+2} \right).
\end{align*}
The first statement of Theorem \ref{thm-eig} is complete by
$$
\begin{aligned}
	&\mathbb{E}(\|\hat\phi_{j}-\phi_{j}\|^2)\lesssim \sum_{k\neq j}\frac{\mathbb{E}\left\{\left(\int\tilde{C}_{0} \phi_{j} \phi_{k}\right)^{2}\right\}}{(\lambda_{j}-\lambda_{k})^2}\\=&\sum_{{k \neq j}\atop{k\leq 2j} }\frac{\mathbb{E}\left\{\left(\int\tilde{C}_{0} \phi_{j} \phi_{k}\right)^{2}\right\}}{(\lambda_{j}-\lambda_{k})^2}+\sum_{k>2j }\frac{\mathbb{E}\left\{\left(\int\tilde{C}_{0} \phi_{j} \phi_{k}\right)^{2}\right\}}{(\lambda_{j}-\lambda_{k})^2}\\
	\lesssim & h^{4}j^{2c+2}+\frac{j^2}{n}\left\{1+\frac{j^{a}}{\bar{N}_2}+\frac{j^{2a}}{\bar{N}^2_2} \right\} + h^{4}j^{2c+1}+\frac{1}{n}\left(j+\frac{j^{a}}{\bar{N}_2h}+\frac{j^{a+1}}{\bar{N}_2}+\frac{j^{2a}}{\bar{N}^2_2h} \right) \\
	\lesssim&\frac{j^{2}}{n}\left\{1+\frac{j^{a}}{\bar{N}_2} +\frac{j^{2a}}{\bar{N}_2^2}  \right\}+\frac{j^{a}}{nh}\left(\frac{1}{\bar{N}_2} + \frac{j^{a}}{\bar{N}_2^2}\right) +h^{4}j^{2c+2},
\end{aligned}
$$
where the second inequality comes from Lemma 7 in \cite{dou2012estimation} and Lemma \ref{lem:Cjk}. The proof for the fixed regular design is similar, and thus we omit the details for brevity.
\end{appendix}

%%%%%%%%%%%%%%%%%%%%%%%%%%%%%%%%%%%%%%%%%%%%%%
%% Support information (funding), if any,   %%
%% should be provided in the                %%
%% Acknowledgements section.                %%
%%%%%%%%%%%%%%%%%%%%%%%%%%%%%%%%%%%%%%%%%%%%%%
% \section*{Acknowledgements}
% The authors would like to thank ...
% 
% The first author was supported by ...
% 
% The second author was supported in part by ...
 
%%%%%%%%%%%%%%%%%%%%%%%%%%%%%%%%%%%%%%%%%%%%%%
%% Supplementary Material, if any, should   %%
%% be provided in {supplement} environment  %%
%% with title and short description.        %%
%%%%%%%%%%%%%%%%%%%%%%%%%%%%%%%%%%%%%%%%%%%%%%
%\begin{supplement}
%\end{supplement}

\bibliographystyle{imsart-nameyear} % Style BST file (imsart-number.bst or imsart-nameyear.bst)
\bibliography{FPC_ref}       

\end{document}